\documentclass[12pt]{amsart}

\usepackage{amsfonts, amsmath, amssymb, amsthm}
\usepackage{mathrsfs}
\usepackage{todonotes}
\usepackage{xcolor}
\usepackage[
	colorlinks,
	linkcolor={black},
	citecolor={black},
	urlcolor={black}
]{hyperref}

\usepackage[margin=1in,marginparwidth=.65in]{geometry}

\usepackage[normalem]{ulem}

\newtheorem{theorem}{Theorem}[section]
\newtheorem{lemma}[theorem]{Lemma}
\newtheorem{corollary}[theorem]{Corollary}
\newtheorem{proposition}[theorem]{Proposition}
\theoremstyle{definition}
\newtheorem{remark}[theorem]{Remark}

  \newtheorem{claim}{Claim}[theorem]
  
\newenvironment{claimproof}[1][Proof of Claim]{\noindent \underline{#1.} }{\hfill$\diamondsuit$}

\newcommand{\bbB}{{\mathbb{B}}}
\newcommand{\bbC}{{\mathbb{C}}}
\newcommand{\bbH}{{\mathbb{H}}}
\newcommand{\bbN}{{\mathbb{N}}}
\newcommand{\bbR}{{\mathbb{R}}}
\newcommand{\bbZ}{{\mathbb{Z}}}

\newcommand{\boldR}{{\mathbf{R}}}

\newcommand{\calD}{{\mathcal{D}}}
\newcommand{\calG}{{\mathcal{G}}}

\newcommand{\calR}{{\mathcal{R}}}
\newcommand{\calS}{{\mathcal{S}}}

\newcommand{\rmd}{{\mathrm{d}}}

\newcommand{\iop}{{\mathbf{i}}}
\DeclareMathOperator{\trace}{tr}
\DeclareMathOperator{\spectrum}{spec}
\DeclareMathOperator{\adj}{adj}

\newcommand{\PT}{{\mathcal{PT}}}
\newcommand{\LP}{{\operatorname{LP}}}
\newcommand{\monodromy}{{\boldsymbol{\Phi}}}
\newcommand{\onestep}{{\mathbf{A}}}
\newcommand{\floqmat}{{\mathbf{H}}}
\newcommand{\idty}{{\mathbf{1}}}
\newcommand{\projector}{{\boldsymbol{\Pi}}}
\DeclareMathOperator{\Var}{Var}
\DeclareMathOperator{\capa}{cap}
\newcommand{\fgset}{{\Sigma}}
\newcommand{\quasimomentum}{{\varkappa}}
\newcommand{\DOSM}{{\kappa}}
\newcommand{\DOS}{{\mathsf{N}}}

\DeclareMathOperator{\real}{Re}
\DeclareMathOperator{\imaginary}{Im}
\renewcommand{\Re}{\real}
\renewcommand{\Im}{\imaginary}

\numberwithin{equation}{section}
\definecolor{cbBlue}{HTML}{0072B2}
\definecolor{cbVermillion}{HTML}{D55E00}
\definecolor{cbGreen}{HTML}{009E73}
\definecolor{cbPurple}{HTML}{CC79A7}

\title[Dispersion for Limit-Periodic Schr\"odinger Operators]{Dispersive Estimates for Limit-Periodic Schr\"odinger Operators}
\author[A. Black \and D. Damanik \and J. Fillman \and G. Young]{Adam Black \and David Damanik \and Jake Fillman \and Giorgio Young}

\address{Department of Mathematics, University of California, Berkeley, Berkeley, CA~94720}
\email{adamblack@berkeley.edu}

\address{Department of Mathematics, Rice University, Houston, TX~77005, USA}
\email{damanik@rice.edu}

\address{Department of Mathematics, Texas A\&M University, College Station, TX~77843, USA}
\email{fillman@tamu.edu}

\address{Department of Mathematics, University of Wisconsin, Madison, WI~53706, USA}
\email{gfyoung@wisc.edu}

\date{}

\newcommand{\vertiii}[1]{{\left\vert\kern-0.25ex\left\vert\kern-0.25ex\left\vert #1 
    \right\vert\kern-0.25ex\right\vert\kern-0.25ex\right\vert}}

\begin{document}

\begin{abstract}
We prove dispersive estimates for limit-periodic operators that are super-exponentially well-approximated by periodic operators. The key ingredient is a sharp enhanced convexity statement for the density of states of periodic Jacobi matrices. This statement implies log-convexity of the density of states and hence is a strengthening of existing results on the convexity of equilibrium measures. We expect that it will be of interest and use outside of this work. We also construct almost-periodic discrete Schr\"odinger operators with purely singular continuous spectrum whose unitary groups enjoy recurrent dispersive motion: on one subsequence of time scales, the unitary group has a dispersive estimate at a rate as close to the free rate as one desires, while on another set of time scales, the unitary group is as far from dispersive as possible. In fact, this behavior is generic within the class of limit-periodic operators.
\end{abstract}

\maketitle

\hypersetup{
	linkcolor={black!30!blue},
	citecolor={black!30!green},
	urlcolor={black!30!blue}
}

\section{Introduction}
\subsection{Main Results}
A basic question in mathematical quantum mechanics is: how do the spectral properties of a Hamiltonian $H$ affect the corresponding quantum dynamics?
The most general result in this direction is the RAGE theorem.
In the context of self-adjoint operators on $\bbZ^{d}$, 
it shows that the states that (in a time-averaged sense) evacuate any compact region under the Schr\"{o}dinger evolution
\begin{align*}
\iop \frac{\partial \psi}{\partial t} = H \psi,
\end{align*}
are precisely those in the continuous spectral subspace of $H$ (see, e.g., \cite[Theorem~5.7]{Teschl2009:SpTh} for a typical statement). Meanwhile, states that are in the absolutely continuous subspace leave compacts even without time-averaging, with this being a one-sided implication rather than a dynamical characterization of spectral type. 
While the RAGE theorem holds in a very general setting, it provides only a qualitative description of the influence of the spectral type on quantum dynamics.
For more precise statements, one is typically interested in quantifying two effects: transport and dispersion.
The former corresponds to the motion of the wavefront, while the latter causes the wave packet to spread. The common expectation, borne out by many rigorous results, is that transport under the flow of $H$ is highly connected to the continuity of its spectral measure, with the extreme case of pure point spectrum suppressing ballistic transport \cite{Simon1990CMP}, and the expectation of ballistic transport in the presence of absolutely continuous spectrum in one space dimension.
Dispersion may be expected when the generalized eigenfunctions have a wave-like character, but it has been most extensively investigated for continuum Schr\"odinger operators with decaying potentials; see \cite{Schlag2007Survey,Schlag2021JMP} for an overview of what is known in this setting, as well as applications to nonlinear PDEs.

Almost-periodic Schr\"odinger operators supply a useful setting for probing this question. 
Indeed, such operators may possess all spectral types, and include the class of periodic operators that are familiar from physical models of crystals; see \cite{DamFil2024ESO2} for relevant background.
In this setting, ballistic transport has been shown by model-dependent means for periodic operators \cite{AschKnauf1998Nonlin, DamLukYes2015CMP, Fillman2021OTAA, BoutetdeMonvelSabri2022OTAA} and for some almost-periodic operators in the absolutely continuous regime \cite{Fillman2017CMP, GeKachkovskiy2023CPAM,  han2025ballistictransportlimitperiodicjacobi, Young2023JST,  Zhao2016CMP, solis2026two}.
In the one-dimensional almost periodic setting, there is also evidence of a direct relationship between absolutely continuous spectrum and ballistic transport without time averaging \cite{ZhangZhao2017CMP}. For more background, see the survey \cite{DamMalYou2026JST} and the books \cite{DamFil2022ESO1, DamFil2024ESO2}.
Less is known about dispersion in this context.
Typically, this effect is quantified through dispersive estimates that measure the decay of the peak of the wave packet. 
Such estimates for continuum periodic (and hence non-decaying) potentials were first proved in \cite{Cai2006DPDE, Cuccagna2008CPDE, Cuccagna2009SIAMJMA, Firsova1987MatSb, Firsova1996JMP}.
It is only recently that dispersive estimates have been established for classes of non-decaying potentials in discrete settings, such as periodic 
\cite{MiZhao2022PAMS, DamFilYou2026FMS, SKW2026JMAA} and quasi-periodic \cite{BambusiZhao2020Advances} potentials.
Other related work in the discrete setting includes dispersive estimates for decaying potentials \cite{EgoKopTes2015JST}, non-locally-finite crystals \cite{KPST2025CMP, LeclercSabriSahlstenPreprint}, trees \cite{AmmariSabri2022AMP}, and a phase-averaged estimate for the Fibonacci Hamiltonian \cite{leclerc2026fourierdecayequilibriumstates}.

In the present work, we consider Jacobi operators $H = H_{a,v}:\ell^2(\bbZ) \to \ell^2(\bbZ)$ given by
\begin{equation}
    [H\psi]_n= a_{n-1} \psi_{n-1} + a_n \psi_{n+1} + v_n \psi_n, \qquad n \in \bbZ,
\end{equation}
where the potential $v:\bbZ \to \bbR$ is bounded and the hopping $a:\bbZ \to (0,\infty)$ is bounded above and below in the sense that $0< \inf a_n \leq \sup a_n < \infty$.
Schr\"odinger operators $H = H_v$ correspond to the case in which $a \equiv 1$.
We denote
\begin{equation}
\vertiii{H} := \max\{\|a\|_\infty, \|a^{-1}\|_\infty, \|v\|_\infty\} < \infty.
\end{equation}

A Jacobi matrix is \emph{periodic} (of period $q \in \bbN$) if $(a_{n+q},v_{n+q}) \equiv (a_n,v_n)$ and \emph{limit-periodic} if $(a,v)$ can be approximated uniformly by periodic sequences (equivalently, $H$ can be approximated by periodic Jacobi matrices in the operator norm topology).  
The Pastur--Tkachenko \cite{PasturTkachenko1,PasturTkachenko2} class $\PT$ consists of those limit-periodic Jacobi matrices $H$ for which there exist $H_k$ of period $q_k$ such that
\begin{equation} \label{eq:PTdef}
    e^{\eta q_{k+1}} \|H - H_k \| \to 0 \quad \text{for all } \eta>0
\end{equation}
 and $q_k | q_{k+1}$ for all $k$.
Such operators are known to have purely absolutely continuous spectral type \cite{Egorova1987} (see also the appendix of \cite{Fillman2017CMP}), Carleson-homogeneous spectrum \cite{FillmanLukic2017JST}, and exhibit a particularly strong form of ballistic transport \cite{Fillman2017CMP}. 
 
We are able to deduce dispersive estimates for operators in $\PT$.
Aside from \cite{BambusiZhao2020Advances}, this is only the second class of examples of almost-periodic operators with aperiodic non-decaying potentials for which such an estimate has been proved. We note that the class in \cite{BambusiZhao2020Advances} is also known to have purely absolutely continuous spectrum.

 \begin{theorem} \label{t:PT}
    If $H \in \PT$, then $H$ obeys a dispersive estimate of the form
    \begin{equation}
        \|e^{-\iop t H}\|_{\ell^1(\bbZ) \to \ell^\infty(\bbZ)}
        \lesssim
        \langle t \rangle^{-\left( \frac13 - \varepsilon \right)}
    \end{equation}
    for every $\varepsilon>0$.
\end{theorem}

Our methods also apply to operators in the exponential class in which one replaces~\eqref{eq:PTdef} by the weaker
\begin{equation} \label{eq:expclassDef}
     \|H - H_k \| \leq A e^{-\eta q_{k+1}}, \text{ for some } A,\eta>0,
\end{equation}
which still implies that $H$ has purely absolutely continuous spectrum if $\eta$ is large enough \cite{solis2026noteacspectrumlimitperiodic}.
In this case, for suitable  $\eta>0$, we obtain dispersion at a rate dependent on $\vertiii{H}$ and $\eta$. See Remark~\ref{rem:EC} for details.

The proof of Theorem~\ref{t:PT} breaks the evolution into time intervals $[T_k, T_{k+1}]$ on which $e^{-\iop t H}$ is very well-approximated by $e^{-\iop t H_k}$ by a Duhamel-type estimate. Naturally, this requires good control on the distance from $H$ to $H_k$, which is supplied by the $\PT$ condition. However, it also requires \emph{uniform} information about the constants in the dispersive estimates from \cite{DamFilYou2026FMS}. More specifically, \cite{DamFilYou2026FMS} guarantees that one has constants $C_1,C_2,\ldots > 0$ such that 
\begin{equation}
    \| e^{-\iop t H_k}\|_{\ell^1(\bbZ) \to \ell^\infty(\bbZ)} \leq C_k \langle t \rangle^{-1/3},
\end{equation}
but without uniform quantitative estimates on $C_k$, one has no information about the $t$ interval on which the decay ``kicks in.'' We emphasize that such estimates are quite delicate and not at all supplied by the arguments of \cite{DamFilYou2026FMS}. 
The difficulty arises because the dispersion constants are controlled by (the minima of) quantities like $|E_j''(k)| + |E_j'''(k)|$, where $E_j$ denotes the $j$th Floquet band function. Thus, one needs quantitative non-degeneracy of the band functions that holds uniformly with respect to all operators $H_k$ with fixed period and bounded $\vertiii{H_k}$, and this is complicated precisely because such band functions can be (and often are) degenerate at the band edges, stymieing the na\"ive perturbative argument. 
Furthermore, this is precisely the phenomenon that one must handle in the current setting, since we are typically perturbing a $q_k$-periodic potential in the space of $q_{k+1}$-periodic potentials, which (typically) involves splitting many degenerate eigenvalues as spectral gaps are opened.

The key ingredient allowing us to produce the necessary control is a new enhanced convexity statement for the density of states of a periodic Jacobi matrix (specifically Theorem~\ref{t:convexity}). 
The new statement is strictly stronger than log-convexity and is exactly what is needed to get precise control over the constants. More precisely, the refined convexity statement enables us to show that $C_k$ can be chosen to be of the form $e^{Cq_k}$ where $C$ is a constant that depends \emph{only} on $\vertiii{H_k}$.
In fact, we prove the indicated convexity statement for the more general class of equilibrium measures (in the sense of potential theory) of finite gap sets, see Theorem~\ref{thm:eqConvexity}.
Convexity in this setting was previously studied by Benko--Dragnev--Totik and Benko--Dragnev \cite{BenDraTot2012RevMatIbe,BenkoDragnev2012ConstrApp}, but their 
techniques only suffice to establish log-convexity.

As a counterpoint to this (and the previously-understood behavior in the periodic case), we show that in the realm of singular continuous spectrum, one can have \emph{recurrent} dispersion in which wave packets disperse at a rate as close to the free rate as one desires on one sequence of time scales while along another sequence the evolved state returns arbitrarily close in modulus to its initial peak.
This result in fact holds for typical limit-periodic Schr\"odinger operators.
To formulate the typicality statement, we denote by 
    $\LP(\bbZ) $  the set of limit-periodic sequences, that is the closure of the set of periodic sequences $\bbZ \to \bbR$ in the uniform topology.
    Since this is a closed subset of the Banach space $\ell^\infty(\bbZ)$, $\LP(\bbZ)$ is itself a complete metric space, so typicality is meant in the sense of the Baire category theorem.
    For additional background about limit-periodic operators, we direct the reader to the survey \cite{DamFil2020LPSurvey}.
    
\begin{theorem} \label{t:LPrecDisp}
For any increasing $\varphi:[0,\infty) \to (0,\infty)$ satisfying
\begin{equation}
    \label{eq:toInfty}
\lim_{t \to \infty} \varphi(t) = \infty,
\end{equation}
 there is a residual set   $\calG = \calG_\varphi \subseteq \LP(\bbZ)$ such that for every $v \in \calG$, the Schr\"odinger operator $H_v$ exhibits purely singular continuous spectrum and $e^{-\iop t H_v}$ enjoys recurrent dispersion in the sense that
\begin{align}
\label{eq:subseqDispersion}
    \liminf_{t \to \infty}  \frac{\langle t \rangle^{1/3}}{\varphi(t)} \|e^{-\iop t H_v}\|_{\ell^1(\bbZ) \to \ell^\infty(\bbZ)} 
    & = 0,\\
    \label{eq:subseqRecurrence}
    \limsup_{t \to \infty}   \|e^{-\iop t H_v} \delta_0 \|_\infty & = 1.
\end{align}

\end{theorem}
Note that \eqref{eq:subseqRecurrence}, in particular, precludes generic limit-periodic Schr\"{o}dinger operators from exhibiting dispersive decay, and in fact, we see that the global estimate fails for the particular state $\delta_0$.
The result above is formulated for generic limit-periodic Schr\"odinger operators, but we expect that an analogous formulation also holds for (generic) Jacobi matrices. Two of the  ingredients needed for our proof of Theorem~\ref{t:LPrecDisp} already exist in the Jacobi setting, namely, dispersion for periodic Jacobi matrices \cite{DamFilYou2026FMS} and generic singular continuous spectrum for Jacobi matrices \cite{DamFilWan2023MathNach, EFGL2022JFA}. The other necessary ingredient is a dense subset of limit-periodic operators exhibiting pure point spectrum. The relevant result for Schr\"odinger operators can be found in \cite{DamGor2016AdvM} and does not exist in the same form for Jacobi matrices. However, since the main ingredient therein is the Kunz--Souillard localization scheme \cite{KunzSouillard1980CMP}, which has been worked out for Jacobi matrices \cite{Bucaj2018OperMat}, we expect that there is a dense subset of limit-periodic Jacobi matrices having pure point spectrum.

Genericity of singular continuous spectrum for limit-periodic operators goes back to work of Avila \cite{Avi2009CMP} and is only stated here for completeness.
Avila's techniques have been extended to other settings, including continuum Schr\"odinger operators \cite{DamFilLuk2017JST} and quantum walks \cite{FillmanOng2017JFA}.

As in the earlier works \cite{BambusiZhao2020Advances, MiZhao2020JMAA,MiZhao2022PAMS}, we may use Theorem~\ref{t:PT} to deduce an $\ell^1 \to \ell^\infty$ dispersive decay estimate for the discrete nonlinear Schr\"odinger (NLS) equation:
\begin{align}\label{eq:dnls}
{\iop}\frac{\partial\psi_n}{\partial t}
=[H\psi]_n\pm |\psi_n|^{\sigma-1}\psi_n,\quad n\in\bbZ,
\end{align}
in the regime of small initial data.
In addition to being a natural discretization of the continuum NLS equation, equation \eqref{eq:dnls} (with periodic potential) has found applications in the physics literature including optical waveguides \cite{eisenberg2002optical,peschel2002optical},
Bose--Einstein condensates \cite{cataliotti2001josephson,cataliotti2003superfluid}, and photonic lattices \cite{efremidis2002discrete,sukhorukov2003spatial}. See also the survey article \cite{kevrekidis2001discrete}. We have the following small-data dispersive estimate, which follows from Theorem~\ref{t:PT} by the same argument that was used to prove \cite[Corollary~2]{BambusiZhao2020Advances}.

\begin{corollary}\label{eq:periodicNLS}
Let $\sigma > 5$ and suppose $H\in \PT$. For every $\varepsilon>0$, there exist constants $C = C(H,\sigma , \varepsilon) , \delta = \delta(H,\sigma ,\varepsilon) > 0$ so that if $\|\psi(0)\|_1\leq \delta$, then $\psi(t)$,  the solution to the Cauchy problem for \eqref{eq:dnls} with initial data $\psi(0)$, satisfies 
\begin{align}\label{eq:nlsdispersion}
\|\psi(t)\|_\infty\leq C\langle t\rangle^{-\left( \frac13 - \varepsilon \right)}.
\end{align}
\end{corollary}

\section{Periodic Jacobi Matrices and a New Convexity Result}

One of the key ingredients in the proof of Theorem~\ref{t:PT} is a quantitative refinement of the main result of \cite{DamFilYou2026FMS} in the periodic setting. This is achieved by way of a convexity statement for {periodic} Jacobi matrices, which we expect to be of interest independently of the rest of the work. The purpose of this section is to recall some background about periodic Jacobi matrices that will be used in the proofs and to formulate this new convexity theorem. We refer the reader to \cite[Section~10.10]{Lukic2022book} for more details and background.

To that end, let us consider a $q$-periodic Jacobi matrix $H = H_{a,v}$.   The one-step transfer matrices, monodromy, and discriminant at spectral parameter $z \in \bbC$ are, respectively, given by
\begin{align}
\onestep(s,t;z) &:= \frac{1}{s} \begin{bmatrix}
    z-t & -1 \\
    s^2 & 0
\end{bmatrix}, \qquad s \neq 0, \ t \in \bbR, \\[2mm]
    \monodromy(z) & := \monodromy_{a,v,q}(z) := \onestep(a_q,v_q;z) \cdots \onestep(a_1,v_1;z), \\[2mm]
    \Delta(z) & := \Delta_{a,v,q}(z) := \trace \monodromy_{a,v,q}(z).
\end{align}
Counted with multiplicity,  $4 - \Delta^2$ has $2q$ real zeros $\lambda_1 \leq \cdots \leq \lambda_{2q}$ and $\Delta'$ has $q-1$ real zeros $x_1 \leq \cdots x_{q-1}$, which are known to obey
\begin{equation} \label{eq:interlace}
    \lambda_{2j-1} < \lambda_{2j} \leq x_j \leq \lambda_{2j+1} < \lambda_{2j+2}
\end{equation}
for every $j=1,2,\ldots,q-1$.
Furthermore, the spectrum of $H$ is given by
\begin{equation}
    \spectrum H_{a,v}
    = \bigcup_{j=1}^q [\lambda_{2j-1}, \lambda_{2j}].
\end{equation}
We call $\sigma_j := [\lambda_{2j-1}, \lambda_{2j}]$ the $j$th \emph{band} of the spectrum.

For $\quasimomentum  \in \bbR$, we denote by $\floqmat(\quasimomentum )$ the Floquet matrix given by restricting $H$ to $\{1,2,\ldots,q\}$ with the boundary condition $\psi_{n+q} = e^{\iop \quasimomentum  q} \psi_n$. For $q \geq 3$, this is given by
\begin{equation}
    \floqmat(\quasimomentum )
    =\begin{bmatrix}
        v_1 & a_1 &&& e^{-\iop \quasimomentum q} a_q \\
        a_1 & v_2 & a_2 \\
        &\ddots & \ddots & \ddots \\
        && a_{q-2} & v_{q-1} & a_{q-1}\\
        e^{\iop \quasimomentum q}a_q &&& a_{q-1} & v_q
    \end{bmatrix}.
\end{equation}
This has eigenvalues (counted with multiplicity) $E_1(\quasimomentum ) \leq \cdots \leq E_q(\quasimomentum )$. The characteristic polynomial of $H(\quasimomentum)$ is related to the discriminant by
\begin{equation} \label{eq:FloquetCharPoly}
    \det(z \idty - H(\quasimomentum))
    = A_q(\Delta(z) - 2 \cos(q\quasimomentum)),
\end{equation}
where $A_q = a_1a_2 \cdots a_q$.
With the help of \eqref{eq:FloquetCharPoly}, we can see that the range of $E_j$ is the $j$th band:
\begin{equation}
    \sigma_j = \{E \in \bbR : E = E_j(\quasimomentum ) \text{ for some } \quasimomentum  \}.
\end{equation}

We then consider the function $ u = u_{a,v,q}$ given by
\begin{equation}  \label{eq:uDef}
    u(z) :=  q^2\frac{4-\Delta(z)^2}{\Delta'(z)^2}.
\end{equation}

In view of \eqref{eq:interlace}, any zero of $\Delta'$ that lies in $ \spectrum H_{a,v}$ is a removable singularity of $u$, so we see that $u$ is meromorphic with poles located at the zeros of $\Delta'$ that lie in the open gaps of $\spectrum H$, and one has the following for each \emph{connected component}  $[\alpha,\beta]$ of $\spectrum H$:
\begin{equation} \label{eq:uOnSpectralComponent}
    u(\alpha) = u(\beta)=0, \qquad u(E)>0, \ E \in (\alpha,\beta).
\end{equation}

\begin{theorem} \label{t:convexity}
    For any $q$-periodic Jacobi matrix $H$ with spectrum $\Sigma$ we have
    \begin{equation} \label{eq:superconvexity}
        u''(E) \leq  -2, \qquad E \in \Sigma^\circ.
    \end{equation}
    This bound is sharp within the class of periodic Jacobi matrices.
\end{theorem}
We will prove the bound \eqref{eq:superconvexity} in Section~\ref{sec:convex}.
Sharpness can be checked by a direct computation: Consider the case  $a_n \equiv a$ and $v_n \equiv v$. One has
\begin{equation}
    \Delta(z) = \frac{z-v}{a}, \qquad u(z) = \frac{4 - a^{-2}(z-v)^2}{a^{-2}} = 4a^2 - (z-v)^2,
\end{equation}
and therefore $u'' \equiv -2$.

 Let us conclude this section by giving a pair of reformulations of Theorem~\ref{t:convexity}.
The function $u$ is closely related to the Marchenko--Ostrovskii map, $\Theta$, defined by \cite[Eq.~(10.47)]{Lukic2022book}, for example. By \cite[Lemma~10.64]{Lukic2022book},
\begin{equation} \label{eq:uToTheta}
    u = [\Theta'(z)]^{-2}.
\end{equation}
Thus, $u'' \leq -2$ is equivalent to $[(\Theta')^{-2}]'' \leq -2$, which can be rearranged to
\begin{equation}
    \Theta'\Theta''' \geq 3 (\Theta'')^2 + (\Theta')^4.
\end{equation}

In a similar vein, let $\DOSM = \DOSM_{a,v}$ denote the {density of states} measure of $H_{a,v}$ (cf.\ \cite[Section~4.3]{DamFil2022ESO1} and \cite[Section~5.4]{Simon2011:Szego}), which is given by
\begin{equation}
    \int \! f \, \rmd \kappa
    = \frac{1}{q} \sum_{j=0}^{q-1} \langle \delta_j, f(H) \delta_j \rangle.
\end{equation}

By standard results (see e.g.\ \cite[Corollary~5.4.20]{Simon2011:Szego}), $\DOSM_{a,v}$ has an analytic density on $\Sigma^\circ$, which we denote by $\DOS =\DOS_{a,v} := \rmd \DOSM / \rmd E$ and is given by
\begin{equation} \label{eq:rhoProdFormula}
    \DOS_{a,v}(E)
    = \frac{1}{\pi} \frac{\prod_{i=1}^{q-1} |E- x_i| }{\prod_{m=1}^{2q} |E - \lambda_m|^{1/2}}.
\end{equation}
Comparing \eqref{eq:rhoProdFormula} and \eqref{eq:uDef}, we see that $\DOS^{-2} = \pi^2 u$, and hence:
\begin{equation}
\label{eq:u''iff}
    u'' \leq -2 \text{ if and only if } (\DOS^{-2})'' \leq -2\pi^2.
\end{equation}

\section{Proof of the Convexity Theorem}
\label{sec:convex}

We will prove Theorem~\ref{t:convexity} by exploiting the connection between the density of states measure of a periodic Jacobi matrix and equilibrium measures in the sense of potential theory. 
We freely use the terminology and notions of potential theory in the complex plane. Further background and proofs can be found in \cite{Saff2010SurvApproxTh} and \cite[Section~5.5]{Simon2011:Szego}.

In this section, we will work in the general setting of finite-gap sets, that is, finite collections of non-degenerate intervals in $\bbR$.
Thus, let $\fgset \subseteq \bbR$ be such a set, and let $\mu_\fgset$ be its equilibrium  measure.
It is well known (see, e.g., \cite[Theorem A.13]{Simon2007IPI}) that $\mu_\fgset$ has a real-analytic density on $\Sigma^\circ$, which we denote by $\rho_\fgset$.
In the setting in which $\fgset = \spectrum H_{a,v}$ is the spectrum of a periodic Jacobi matrix, it is known (\cite[Theorem~5.5.17]{Simon2011:Szego}) that $\mu_\fgset = \DOSM_{a,v}$, the density of states measure of $H_{a,v}$, so $\rho_\fgset = \DOS_{a,v}$, the density of states.

In this generality, we will show that $
\rho_\fgset$ satisfies $(\rho_\fgset^{-2})''\leq -2\pi^2$ on $\Sigma^\circ$.
Since not all finite-gap sets arise as spectra of periodic Jacobi matrices (see \cite[Theorem~5.5.25]{Simon2011:Szego}), the following theorem is strictly stronger than needed for Theorem~\ref{t:convexity}:

\begin{theorem}
\label{thm:eqConvexity}
Let $\Sigma$ be a finite-gap set and $\rho_\Sigma$ the density of the equilibrium measure of $\Sigma$ as above.
For every $x\in \fgset^{\circ}$, we have that
\begin{align}
\label{eq:convexityBound}
\frac{\rmd^2}{\rmd x^2}\bigl(\rho_\fgset(x)^{-2}\bigr)
\leq -2\pi^2.
\end{align}
\end{theorem}

\begin{remark}
    Consider an open interval $I \subseteq \bbR$ and a $C^2$ function $f:I \to (0,\infty)$, and notice that
    \(
        (f^{-2})'' \leq -2\pi^2
    \)
    implies
    \begin{equation}
        f''f-3(f')^2 \geq \pi^2 f^4,
    \end{equation}
    which in turn implies the following $\log$-convexity statement:
    \begin{equation}
        (\log f)''
        = \frac{f''f-(f')^2}{f^2}
        \geq \frac{2(f')^2 + \pi^2 f^4}{f^2}
        = 2 [(\log f)']^2 + \pi^2 f^2
        >0.
    \end{equation}
    However, $\log$-convexity of $f$ alone does not guarantee $(f^{-2})'' \leq -2\pi^2$ or even $(f^{-2})'' < 0$ (e.g., consider $f(x) = \exp(ax^2)$ with $a>0$).
    Thus, \eqref{eq:convexityBound} is a strictly stronger statement than $\log$-convexity of $\rho_\Sigma$, which was established in \cite[Theorem~1.1]{BenDraTot2012RevMatIbe}.
 More generally,
inequalities of the form \eqref{eq:convexityBound} have been studied in the elliptic and parabolic PDE literature under the moniker \emph{power convexity/concavity}, see \cite{Kennington1985} for a classical reference on this topic and \cite{BianchiniSalani2013} for more modern developments.
 If $f: U\to (0,\infty)$ is a positive $C^2$ function on a domain $U \subset \bbR^d$ and $\alpha >0$, then $f$ is said to be \emph{$\alpha$-convex} if $f^\alpha$ is convex and $0$-convex if $f$ is $\log$-convex. When $\alpha<0$, $f$ is said to be $\alpha$-convex if $f^\alpha$ is concave.
The smaller the $\alpha$, the stronger the statement, as can be seen from the equivalence
\begin{align*}
f\text{ is $\alpha$-convex}
\iff
D_v^2f \cdot f + (\alpha -1) (D_v f)^2 \geq 0, \quad \text{for all } v\neq 0,
\end{align*}
where $D_v$ denotes the directional derivative in the direction $v$ (cf.\ \cite[Property~4]{Kennington1985}, noting that \cite{Kennington1985} discusses $\alpha$-concavity and we are concerned with $\alpha$-convexity).
The inequality \eqref{eq:convexityBound} therefore implies that  $\rho_{\Sigma}$ is strictly $-2$-convex.
Due to the example \eqref{eq:sqrtDensity} below, this $-2$ is actually the sharp $\alpha$ for which $\alpha$-convexity can be proven for equilibrium measures of intervals. 
\end{remark}

Notice that, if $f > 0$ is a $C^2$ function on an interval $I$ such that $\frac{\rmd^2}{\rmd x^2}\bigl(f(x)^{-2}\bigr) = -2\pi^2$ holds on $I$, then $f$ must be of the form 
\begin{align}
\label{eq:sqrtDensity}
f(x) =  \frac{1}{\pi \sqrt{(x-a)(b-x)} }, \quad x \in I,
\end{align}
which is the density of the equilibrium measure of the interval $[a,b]$ (see \cite[Example~1.11]{Saff2010SurvApproxTh}).
Motivated by this, the idea will be to compare $\rho_\fgset$ to such a density.

We will use the relationship between $\rho_\fgset$ and the Green's function of $\fgset$ from infinity.
To this end, define the Green's function $g_\fgset$ by
\begin{align*}
g_\fgset(z) = -\log(\capa(\fgset)) - \int_\fgset \log \frac{1}{|z-y|} \: \rmd \mu_\fgset(y), \quad z \in \bbC.
\end{align*}
We summarize its relevant properties, all of which are standard.

For a function $g: \bbC \to \bbR$, we view $z = x+\iop y$ and identify $g(z)$ with  $g(x,y)$, a function of two real variables. As usual, expressions like $g(x,0^+)$ and $\partial_y g(x,0^+)$ denote limits, for example,
\begin{equation}
    \partial_y g(x,0^+) = \lim_{\varepsilon \downarrow 0} \partial_y g(x,\varepsilon).
\end{equation}

\begin{proposition}
\label{pr:eqFacts}
The Green's function $g_\fgset$ has the following properties:
\begin{enumerate}
\item The function $g_\fgset$ is continuous on all of $\bbC$, is harmonic and non-negative on $\bbC\setminus \fgset$, and vanishes on $\fgset$.
As $z\to \infty$,
\begin{align}
\label{eq:gAsy}
g_\fgset(z) = \log|z| - \log \capa(\fgset) + O(|z|^{-1}).
\end{align}
Moreover, it satisfies $g_\fgset(z) = g_\fgset(\overline{z}) $.
\item For every $x\in \fgset^{\circ}$, we have that
\begin{align*}
\partial_y g_\fgset(x,0^{+}) = \pi \rho_\fgset(x).
\end{align*}
\end{enumerate}
\end{proposition}
\begin{proof}
Frostman's theorem \cite[Theorem A.IV.2.1]{SaffTotik2017} shows that $g_\fgset$ is non-negative on $\bbC$ and equal to $0$ quasi-everywhere on $\fgset$.
The continuity of $g_\fgset$ on $\bbC\setminus \fgset$ follows directly from the integral representation, whereas on $\Sigma$ it follows from the regularity of the Dirichlet problem for finite-gap sets, see the discussion after \cite[Eq.~(A.42)]{Simon2007IPI}.
Continuity then implies that $g_\fgset=0$ \emph{everywhere} on $\fgset$. 
The fact that $g_\fgset$ is harmonic on $\bbC\setminus \fgset$, the asymptotic \eqref{eq:gAsy}, and the symmetry $g_\fgset(z) =g_\fgset(\overline{z})$ all follow easily from the integral representation.
The identity $\partial_y g_\fgset(x,0^{+}) = \pi \rho_\fgset(x)$ is   \cite[Eq.~(A.51)]{Simon2007IPI}.
\end{proof}

The proof of Theorem~\ref{thm:eqConvexity} will proceed as follows.
The inequality \eqref{eq:convexityBound} is equivalent to the fact that on each connected component of $\fgset$, $\rho^{-2}_\fgset(x) + \pi^{2}x^{2}$ lies above all of its secant lines.
Given $r<s$ fixed in the interior of $\fgset$, we will choose an interval $I$ such that the density of its equilibrium measure, $\rho_I$, coincides with $\rho_\fgset$ at $r$ and $s$.
Since $ \frac{\rmd^2}{\rmd x^2}\bigl(\rho_I(x)^{-2}\bigr) = -2\pi^{2}$ (and therefore $\rho_I^{-2} + \pi^{2}x^{2}$ is equal to the secant line), the claim will follow if we can show that $\rho_\fgset \leq \rho_I$ on $[r,s]$.
This is a consequence of the following more general comparison result for harmonic functions:
\begin{proposition}
\label{pr:harmonic}
Let $\fgset$ be a finite-gap set and $I \subset \bbR$ be a compact interval such that $I \cap \Sigma^\circ \neq \emptyset$,  and let $D:\hat{\bbC} \to \bbR$ be harmonic on $\hat{\bbC}\setminus \left( \fgset \cup I \right) $, be continuous on $\hat{\bbC}$, and satisfy the symmetry condition $D(z) = D(\overline{z})$.
Fix $[r,s] \subset I \cap \fgset^{\circ}$, and suppose that $D$ satisfies
\begin{itemize}
\item $D\geq 0$	on $I$.
\item  $D\leq 0$ on $\fgset$.
\item $\partial_y D(r,0^{+})\leq 0$ and $\partial_y D(s,0^{+})\leq 0$.
\end{itemize}
Then $\partial_y D(x,0^{+}) \leq 0$ for all $x\in [r,s]$.
\end{proposition}

\begin{proof}
Suppose that $\partial_y D(x_0,0^{+}) >0$ for some $x_{0} \in [r,s]$,  so that $D(x_0,y)>0$ for all $y$ sufficiently small.
Let $U$ be the connected component of 
\begin{align*}
\{z \in \bbH  \mid  D(z) >0\},
\end{align*}
containing a segment above $x_0$.

Suppose that $D$ attains its maximum on $\overline{U}$ at the point $\hat{z}$.
We claim that $\hat z \in I \setminus \Sigma$.
By the strong maximum principle, $\hat{z}$ must lie on the boundary of $U$.
If $\hat{z}$ occurred in $\bbH$, then we would have that $D$ is positive on some neighborhood of $\hat{z}$ contained in $\bbH$, and therefore $\hat{z}$ would not be on the boundary of $U$.
We must then have that $\hat{z} \in \bbR$.
If $\hat{z} \in  \bbR\setminus(\fgset\cup I)$ then by symmetry we may find a disc centered at $\hat{z}$ and contained in $\bbC\setminus (\fgset\cup I)$, on which $D$ is maximized at $\hat{z}$.
But this contradicts the maximum principle, so we must have that $\hat{z} \in  \fgset \cup I$.
Now, if $\hat{z} \in \fgset$, then $D(\hat{z})\leq 0$, which contradicts the fact that $U$ is non-empty. 
Therefore, we must have that $\hat{z} \in I \setminus \fgset\subset I \setminus [r,s]$, so suppose without loss of generality that $\hat{z} <r$.

Now choose a smooth simple path $\gamma$ connecting $x_0$ to $\hat{z}$ entirely through $U$. 
Let $\Omega$ be the domain with boundary given by $\gamma$ and the segment $[\hat{z},x_0]$.
See Fig.~\ref{fig:harmonic} for an illustration.

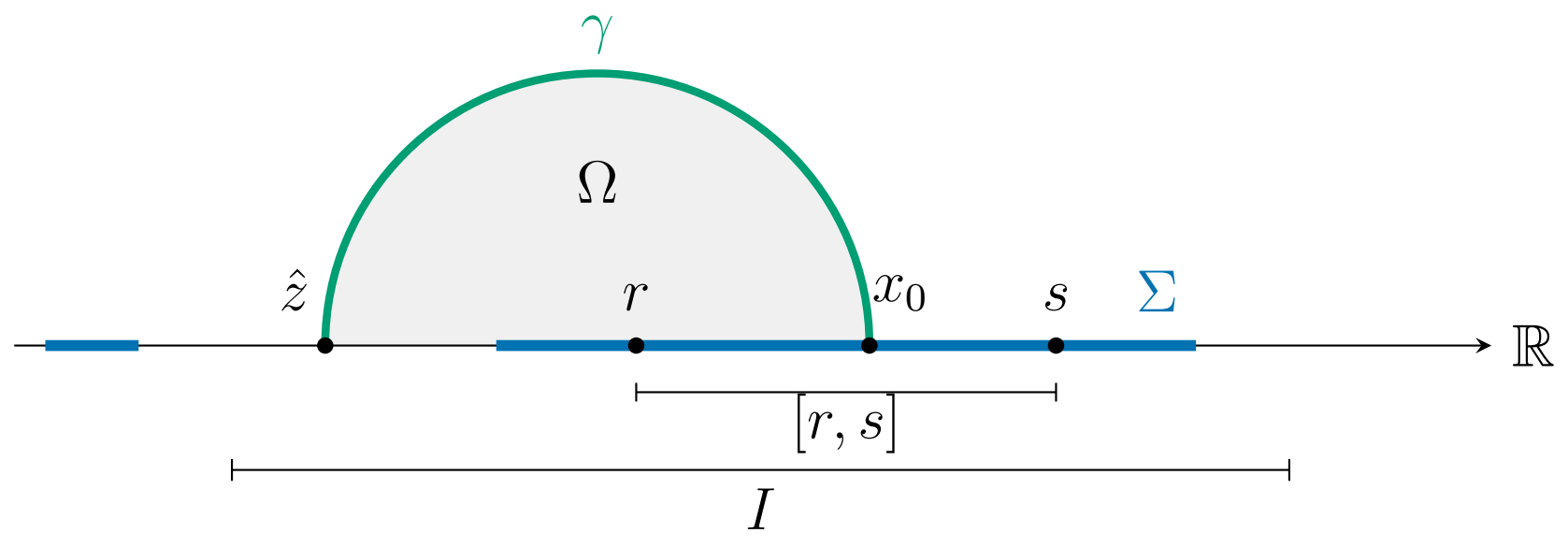
\begin{figure}[h]
\begin{tikzpicture}[scale=1.1, >=stealth]

  \fill[gray!12]
    (-3,0) arc[start angle=180,end angle=0,radius=1.75]
    -- cycle;
  \draw[ultra thick, color = cbGreen]
    (-3,0) arc[
    start angle=180,end angle=0,radius=1.75];
  \node at (-1.25,1.05) {$\Omega$};
  \node at (-1.25,2.00) {\textcolor{cbGreen}{$\gamma$}};

  \draw[->] (-5,0) -- (4.5,0) node[right] {$\bbR$};

  \draw[line width=2.2pt,color=cbBlue] (-4.8,0) -- (-4.2,0);
  \draw[line width=2.2pt,color=cbBlue] (-1.9,0) -- (2.6,0);
  \node at (2.35,0.35) {\textcolor{cbBlue}{$\fgset$}};

  \foreach \x/\name in {
    -3/\hat{z},
    -1/r,
    0.5/x_0,
    1.7/s
  } {
    \fill (\x,0) circle (1.5pt);
  }
\node[above=3pt] at (-3.2,0){$\hat{z}$};
\node[above=3pt] at (-1,0){$r$};
\node[above=3pt] at (0.7,0){$x_0$};
\node[above=3pt] at (1.7,0){$s$};

  \draw[thin] (-1,-0.3) -- (1.7,-0.3);
  \draw (-1,-0.24) -- (-1,-0.36);
  \draw (1.7,-0.24) -- (1.7,-0.36);
  \node at (0.35,-0.5) {$[r,s]$};

  \draw[thin] (-3.6,-0.8) -- (3.2,-0.8);
  \draw (-3.6,-0.73) -- (-3.6,-0.87);
  \draw (3.2,-0.73) -- (3.2,-0.87);
  \node at (-0.2,-1.05) {$I$};

\end{tikzpicture}
\caption{Illustration of some points and sets in the proof of Proposition~\ref{pr:harmonic}.}
\label{fig:harmonic}
\end{figure}

We have that $[\hat{z},x_0] \subset I$ so $D\geq 0$ there, and therefore $D\geq 0$ on $\partial \Omega$.
By the strong minimum principle, $D >0$ on $\Omega$.
But then, $D(z) > 0$ for all $z$ in $\Omega$ sufficiently close to $r$, which by the Hopf lemma (see, e.g., \cite[Lemma~14.5.1]{GreeneKrantz2002GSM}) implies that $\partial_y D(r,0^{+}) > 0$.
This contradicts the assumption that  $\partial_y D(r,0^{+}) \leq 0$, so we are done.
\end{proof}

\begin{proof}[Proof of Theorem~\ref{thm:eqConvexity}]
Fix $[r,s] \in \fgset^{\circ}$.
Choose $a$ and $b$ such that 
\begin{align*}
\rho_I := \frac{1}{\pi \sqrt{(x-a)(b-x)} }
\end{align*}
satisfies $\rho_\fgset(r) = \rho_I(r)$ and $\rho_\fgset(s) = \rho_I(s)$.
Note that $\rho_I$ is the density of the equilibrium measure of the interval $I = [a,b]$.
Per the above discussion, it suffices to show that $\rho_\fgset \leq \rho_I$ on $[r,s]$.

Let $g_I$ be the Green's function of $I$, and let $D(z) = g_\fgset(z) - g_I(z)$.
From Proposition~\ref{pr:eqFacts}, it is clear that $D$ satisfies all of the conditions of Proposition~\ref{pr:harmonic}. 
In particular it is harmonic at infinity because of \eqref{eq:gAsy}.
We conclude that $\partial_y (g_\fgset(x,0^{+})- g_I(x,0^{+})) \leq 0$, for all $x\in [r,s]$, which again per Proposition~\ref{pr:eqFacts} yields $\rho_\fgset \leq \rho_I$ on $[r,s]$, and this is exactly what we wanted to show. 
\end{proof}

\begin{proof}[Proof of Theorem~\ref{t:convexity}]
From \eqref{eq:u''iff}, it suffices to show that $(\rho_\fgset^{-2})''\leq -2\pi^2$ on $\fgset^\circ$.
However, this is just the special case of Theorem~\ref{thm:eqConvexity} when $\Sigma = \spectrum H$, so we are done.
\end{proof}

\section{Improved Dispersive Bounds for Periodic Operators}

Consider a $q$-periodic Jacobi matrix $H= H_{a,v}$ with band functions $E_j(\cdot)$ and bands $\sigma_j$.
Our main aim in the present section is to prove a quantitative nondegeneracy result on the derivatives of the Floquet functions.
As usual, when $E_j$ is degenerate at $\quasimomentum_\star = 0$ or $\quasimomentum_\star = \pi/q$, we write $E_j^{(r)}({\quasimomentum}_\star)$ for the derivative of the analytic continuation of $E_j$ through ${\quasimomentum}_\star$.
Following \cite{DamFilYou2026FMS}, we fix $q \in \bbN$ and define for $(a,v):\bbZ \to (0,\infty) \times \bbR$ of period $q$ the function
\begin{equation}
    \delta(a,v) = \min_{ 1 \le j \le q} \min_{{\quasimomentum} \in [0,\pi/q]} \big(|E_j''({\quasimomentum})| + |E_j'''({\quasimomentum})| \big) .
\end{equation}
The main technical result of \cite{DamFilYou2026FMS} is the nondegeneracy statement $\delta(a,v) >0$ for all periodic $(a,v)$.
Our next result gives a quantitative refinement:

\begin{theorem} \label{t:deltaLB}
Assume $(a,v)$ is $q$-periodic, let $1 \le j \le q$ be given, and let $[\alpha, \beta] \supseteq \sigma_j$ denote the connected component of $\spectrum H_{a,v}$ containing $\sigma_j$. 
For all ${\quasimomentum} \in [0,\pi/q]$, we have
\begin{equation} \label{eq:Ej''Ej'''LB}
    |E_j''({\quasimomentum})| + |E_j'''({\quasimomentum})| 
    \geq \frac{\beta - \alpha}{2}
    \geq \frac{1}{2} |\sigma_j|.
\end{equation}
In particular,
\begin{equation}
    \delta(a,v) \geq \frac{1}{2} \min_j |\sigma_j|.
\end{equation}
\end{theorem}

Note that for each $j$, the left-hand side of \eqref{eq:Ej''Ej'''LB} is continuous as a function of ${\quasimomentum} \in [0,\pi/q]$, so it suffices to prove \eqref{eq:Ej''Ej'''LB} for ${\quasimomentum} \in (0,\pi/q)$.

We begin by rewriting the derivatives of the band functions in terms of $u$.

\begin{lemma} \label{lem:Ej(r)u}
For all ${\quasimomentum} \in (0,\pi/q)$, we have
\begin{align}
\label{eq:Ej(1)u}
E_j'({\quasimomentum})^2 & = u( E_j({\quasimomentum})),\\
\label{eq:Ej(2)u}
E_j''({\quasimomentum}) & = \frac{1}{2}u'( E_j({\quasimomentum})), \\ 
\label{eq:Ej(3)u}
    E_j'''({\quasimomentum}) & = \frac{1}{2}u''( E_j({\quasimomentum})) E_j'({\quasimomentum}).
\end{align}
In particular,
\begin{equation} \label{eq:E'''GeqE'}
    |E_j'''({\quasimomentum})| \geq |E_j'({\quasimomentum})|.
\end{equation}
\end{lemma}

\begin{proof}
Beginning with $\Delta(E_j({\quasimomentum})) = 2 \cos(q{\quasimomentum})$ and differentiating in ${\quasimomentum}$ gives
\begin{equation}
    \Delta'(E_j({\quasimomentum})) E_j'({\quasimomentum}) = -2q \sin(q{\quasimomentum}).
\end{equation}
Squaring, rearranging, and using the definition of $u$ yields
\begin{equation}
    E_j'({\quasimomentum})^2 = \frac{4q^2 \sin^2(q{\quasimomentum}) }{\Delta'( E_j({\quasimomentum}))^2} = \frac{q^2(4 - \Delta( E_j({\quasimomentum}))^2)}{\Delta'( E_j({\quasimomentum}))^2} = u( E_j({\quasimomentum})),
\end{equation}
as promised. Differentiating \eqref{eq:Ej(1)u} with respect to ${\quasimomentum}$ two times gives \eqref{eq:Ej(2)u} and \eqref{eq:Ej(3)u}.

Since $u'' \leq -2$ by Theorem~\ref{t:convexity}, \eqref{eq:Ej(3)u} yields
\begin{equation}
    \frac{|E_j'''({\quasimomentum})|}{|E_j'({\quasimomentum})|}
    = \frac{\frac{1}{2} |u''(E)||E_j'({\quasimomentum})|}{|E_j'({\quasimomentum})|}
    \geq 1,
\end{equation}
which proves \eqref{eq:E'''GeqE'}.
\end{proof}

Together with the convexity statement $u'' \leq -2$, we now have everything needed to prove the desired quantitative nondegeneracy result.

\begin{proof}[Proof of Theorem~\ref{t:deltaLB}]
 Fix $(a,v)$ of period $q$, and let $\Sigma = \spectrum H_{a,v}$.
 \begin{claim}
On $[\alpha,\beta]$, $u$ has a nondegenerate maximum value satisfying
 \begin{equation}
     u(E_*) \geq \frac{(\beta-\alpha)^2}{4}.
 \end{equation}
 \end{claim}

 \begin{claimproof}
Recalling \eqref{eq:uOnSpectralComponent}, we have $u > 0$ on $(\alpha,\beta)$, $u(\alpha) = u(\beta) = 0$.
Thus, the existence of a nondegenerate maximum at some $E_* \in (\alpha,\beta)$ follows from $u'' \leq-2$.

Furthermore, if we consider
\[ \widetilde u(E) = u(E) + (E-\alpha)(E-\beta), \]
we have $\widetilde{u}(\alpha) = \widetilde{u}(\beta)=0$ and $\widetilde u'' \leq 0$ on $[\alpha,\beta]$.
Consequently, $\widetilde{u} \geq 0$ throughout $[\alpha,\beta]$ (which can be seen, e.g., with the mean value theorem). From this, we deduce
\begin{align}\nonumber
    u(E_*) 
    \geq u(\tfrac{\alpha + \beta}{2})
    & = \widetilde{u}(\tfrac{\alpha + \beta}{2}) - (\tfrac{\alpha + \beta}{2}-\alpha)(  \tfrac{\alpha + \beta}{2} - \beta) \\
    \nonumber
    & \geq (\tfrac{\alpha + \beta}{2}-\alpha)( \beta - \tfrac{\alpha + \beta}{2}) \\
    & =  \frac{(\beta-\alpha)^2}{4},
\end{align}
as desired.
 \end{claimproof}
\bigskip

In view of Lemma~\ref{lem:Ej(r)u}, we have (abbreviating $E = E_j({\quasimomentum})$):
\begin{equation}
    |E_j'''({\quasimomentum})|^2 + |E_j''({\quasimomentum})|^2
    \geq |E_j'({\quasimomentum})|^2 + |E_j''({\quasimomentum})|^2
    = u(E) + \frac{1}{4}[u'(E)]^2,
\end{equation}
and we may note that the derivative of the final expression is given by
\begin{equation}
    \frac{\rmd}{\rmd E}\left(u + \tfrac14 (u')^2 \right)
    = u' + \tfrac{1}{2}u'u'' = \left( 1 + \tfrac{1}{2}u'' \right)u'.
\end{equation}
Since $u$ has a single nondegenerate maximum at $E_*$, we have $u'>0$ on $(\alpha,E_*)$ and $u'<0$ on $(E_*,\beta)$, while Theorem~\ref{t:convexity} implies $1+u''/2\leq 0$.
Putting those two pieces of information together, we therefore see that $u+(u')^2/4$ is \emph{minimized} at $E_*$, leading to
\begin{equation}
    |E_j'''({\quasimomentum})|^2 + |E_j''({\quasimomentum})|^2
    \geq u(E_*) + \frac{1}{4}[u'(E_*)]^2
    \geq \frac{(\beta - \alpha)^2}{4} + 0.
\end{equation}
By the Cauchy--Schwarz inequality, we get
\begin{equation}
    |E_j'''({\quasimomentum})| + |E_j''({\quasimomentum})|
    \geq \left(u(E_*) + \frac{1}{4}[u'(E_*)]^2\right)^\frac{1}{2}
    \geq \frac{\beta - \alpha}{2}
    \geq \frac{1}{2} |\sigma_j|,
\end{equation}
as desired.
\end{proof}

This can now be used to prove a quantitative dispersive estimate with a constant of the form $e^{Cq}$ where $C$ depends \emph{only} on $\vertiii{H}$.
\begin{theorem} \label{t:quant:dispersion}
    Given $R>0$, there exists $C = C(R)>0$ such that every $q$-periodic Jacobi matrix $H$ with $\vertiii{H}\leq R$
    obeys a dispersive estimate of the form
    \begin{equation}
        \|e^{-\iop t H_{a,v}} \|_{\ell^1(\bbZ) \to \ell^\infty(\bbZ)} \leq e^{Cq} \langle t \rangle^{-1/3}.
    \end{equation}
\end{theorem}

\begin{proof}
Fix $H = H_{a,v}$ of period $q$ satisfying $\vertiii{H} \leq R$.
    The basic analysis is essentially the same as the proof of \cite[Theorem~1.1]{DamFilYou2026FMS}, with improvements facilitated by the new non-degeneracy result and an improved estimate on the variation of the amplitude.
    First, recall \cite[Lemma~3.1]{DamFilYou2026FMS} and write
    \begin{equation}
        \langle  e^{-\iop t H} \psi , \delta_{m + \ell q} \rangle 
        = \sum_{m'=1}^q \sum_{\ell' \in \bbZ} \psi_{m' + \ell' q} \int_{\bbB} \sum_{j=1}^q e^{-\iop (t E_j({\quasimomentum}) + q(\ell'- \ell)\quasimomentum)} \Pi_{j}^{m,m'}({\quasimomentum}) \, \frac{\rmd {\quasimomentum}}{|\bbB|},
    \end{equation}
    where $\Pi_{j}^{m,m'}({\quasimomentum}) = {[v_j({\quasimomentum})]_m} \overline{[v_j({\quasimomentum})]_{m'}}$ and $\bbB = [-\pi/q,\pi/q)$.

Combining Theorem~\ref{t:deltaLB} with Theorem~\ref{t:bandLB}, there are constants $B,B'>0$ depending only on $R$ such that 
\begin{equation} \label{eq:bandLB}
    \delta(a,v) \geq \min_j \frac{1}{2} |\sigma_j| \geq \frac{1}{2e^{Bq}} \geq \frac{2}{e^{B'q}}\qquad \text{for all }j.
\end{equation}

Fix $j \in \{1,2,\ldots, q\}$. On account of Lemma~\ref{lem:Ej(r)u}, we see that $E_j'''({\quasimomentum})$ does not vanish in $(0,\pi/q)$ (in particular because $E_j'$ does not vanish), so $E_j''({\quasimomentum})$ is monotonic on $(0,\pi/q)$, and the set of ${\quasimomentum}$ for which $|E_j''({\quasimomentum})| < \delta(a,v)/2$ is an interval (on which one has $|E_j'''({\quasimomentum})| \geq \delta(a,v)/2$ by the definition of $\delta(a,v)$).
Thus, we can partition $[0,\pi/q]$ into (no more than) three subintervals on which either $|E_j''| \geq \delta/2 \geq e^{-B'q}$ or $|E_j'''| \geq \delta/2 \geq e^{-B'q}$.

To conclude, apply the van der Corput Lemma on each of these three intervals and use the following claim:

\begin{claim}
    $\|\Pi_{j}^{m,m'}\|_\infty \leq 1$ and $\operatorname{Var}(\Pi_{j}^{m,m'},[0,\pi/q]) \leq Cq$ for a universal constant $C>0$.
\end{claim}

\begin{claimproof}
    For ${\quasimomentum} \in (0,\pi/q)$,  $\floqmat({\quasimomentum})$ has simple spectrum, so $\Pi_{j}^{m,m'}({\quasimomentum})$ is the $m,m'$ matrix element of $\projector_j({\quasimomentum})$, the spectral projector onto $\ker( E_j(\quasimomentum)\idty - \floqmat(\quasimomentum) )$. 
    Since $\|\projector_j({\quasimomentum})\|=1$, the first estimate follows directly.
 Furthermore, since $|\Pi_{j}^{m,m'}({\quasimomentum}) | \leq 1$ for all $j,{\quasimomentum},m,m'$, to prove the second statement, it suffices to show that $\Re \Pi_{j}^{m,m'}$ and $\Im \Pi_{j}^{m,m'}$ have at most $Cq$ intervals of monotonicity (the argument below gives $C = 3$ for the real part and $C=4$ for the imaginary part).
    
    For ${\quasimomentum}\in (0,\pi/q)$, we notice that \eqref{eq:FloquetCharPoly} gives
    \begin{equation}
        \projector_j({\quasimomentum})  = \lim_{z \to E_j({\quasimomentum})} (z-E_j({\quasimomentum})) (z \idty -\floqmat({\quasimomentum}))^{-1}
        = \frac{\adj(E_j({\quasimomentum})\idty - \floqmat({\quasimomentum}))}{A_q \Delta'(E_j({\quasimomentum}))},
    \end{equation}
    where $\adj(A)$ denotes the adjugate matrix of $A$ (i.e., the transpose of the matrix of cofactors).

The entries of $\boldR_j({\quasimomentum}) := \adj(E_j({\quasimomentum}) \idty  - \floqmat({\quasimomentum}))$ all have the form 
\begin{align*}
     [\boldR_j({\quasimomentum})]_{m,m'}
    & = F_1^{m,m'}(E_j({\quasimomentum})) + F_2^{m,m'}(E_j({\quasimomentum})) \cos q{\quasimomentum}  \\
    & \qquad +\iop F_3^{m,m'}(E_j({\quasimomentum})) \sin q{\quasimomentum} \\
    & =
    F_1^{m,m'}(E_j({\quasimomentum})) 
    + F_2^{m,m'}(E_j({\quasimomentum})) \frac{\Delta(E_j({\quasimomentum}))}{2} \\
    &  \qquad 
    + \iop F_3^{m,m'}(E_j({\quasimomentum})) \frac{\sqrt{4 - \Delta(E_j({\quasimomentum}))^2}}{2},
\end{align*}
where each $F_i^{m,m'}$ is a polynomial of degree $q-1$ or less with real coefficients.
Thus, we can change variables and view $\Pi_j^{m,m'}$ as a function of $E \in \sigma_j$ rather than ${\quasimomentum} \in [0,\pi/q]$.
Notice that $E_j({\quasimomentum})$ is monotonic on $[0,\pi/q]$, so the real and imaginary parts of $\Pi_{j}^{m,m'}({\quasimomentum})$ will have the same number of monotonicity intervals as the real and imaginary parts of the function
\begin{equation}
    \widetilde \Pi_{j}^{m,m'}(E)
    = \frac{F_1^{m,m'}(E)
    +F_2^{m,m'}(E)\frac{\Delta(E)}{2} + \iop F_3^{m,m'}(E) \frac{\sqrt{4-\Delta(E)^2}}{2}}{\Delta'(E)}.
\end{equation}
Rewrite this as
\begin{equation}
    \widetilde{\Pi}_j^{m,m'}(E) = \frac{P_1^{m,m'}(E)  +  \iop P_2^{m,m'}(E)\sqrt{Q(E)}}{\Delta'(E)}
\end{equation}
where $P_1$, $P_2$,   and $Q$ are polynomials of degree at most  $2q-1$, $q-1$, and $2q$, respectively.

 Dropping $m$, $m'$, and $E$ from the notation, we have
\begin{align}
    [\Re \widetilde{\Pi}_j]' 
    & = \frac{P_1' \Delta' - P_1 \Delta''}{[\Delta']^2},\\
    [\Im \widetilde{\Pi}_j]' 
    & = \frac{(P_2'Q^{1/2} + \frac12 P_2 Q^{-1/2}Q')\Delta' - P_2Q^{1/2}\Delta''}{[\Delta']^2}
    = \frac{(P_2'Q + \frac12 P_2 Q' )\Delta' - P_2Q\Delta''}{[\Delta']^2 Q^{1/2}}.
\end{align}
If the numerator of either vanishes identically, the corresponding function is constant. 
The numerator of $[\Re \widetilde{\Pi}_j]'$ is a polynomial of degree at most $3q-3$ and hence has no more than $3q-3$ zeros if it does not vanish identically.
It follows that $\Re \widetilde{\Pi}_j$ has no more than $3q-2$ monotonicity intervals, so $\Var [\Re \widetilde{\Pi}_j] \leq 6q$ (since $-1 \leq \Re \widetilde{\Pi}_j \leq 1$).
In a similar way, we see that $\Var[\Im \widetilde{\Pi}_j] \leq 8q$.
\end{claimproof}
\medskip

As discussed above, applying the van der Corput Lemma  \cite[Ch.~VIII, §1.2]{Stein1993HarmonicAnalysis} on the three intervals constructed above and applying the claim produces the desired bound.
\end{proof}

As a consequence of the sharpened estimates above, we can deduce local boundedness of the optimal dispersion constant on Jacobi matrices of a fixed period. More precisely, for
  a given periodic $(a,v)$ (which we may identify with an element of $(0,\infty)^p \times \bbR^p$ in a canonical manner), let 
    \begin{equation}
        M(a,v) = \sup_{t} \langle t \rangle^{1/3} \|e^{-\iop t  H_{a,v}}\|_{\ell^1(\bbZ) \to \ell^\infty(\bbZ)},
    \end{equation}
which is the minimal $M$ for which 
\[
\|e^{-\iop t H_{a,v}}\|_{1\to\infty} \leq M\langle t \rangle^{-1/3}\]
holds for all $t$.

\begin{corollary} \label{coro:dispConstLocallyBounded}
   For each $p \in \bbN$, the function $M:(0,\infty)^p \times \bbR^p \to (0,\infty)$ is locally bounded.
   That is,
   \begin{equation}
       \sup_{(a,v)\in K} M(a,v) < \infty
   \end{equation}
   for any compact set $K \subseteq (0,\infty)^p \times \bbR^p$.
\end{corollary}

\begin{proof}
This follows immediately from Theorem~\ref{t:quant:dispersion}.
\end{proof}

\section{Dispersion for Limit-Periodic Operators}
 
We are now ready to give the proofs of the main results. In the first proof, we make use of the following simple Duhamel-type estimate for bounded self-adjoint $A$ and $B$ and $\psi \in \ell^1(\bbZ)$:
\begin{equation} \label{eq:duhamel}
    \|(e^{-\iop t A} - e^{-\iop t B}) \psi\|_\infty
    \leq  \|(e^{-\iop t A} - e^{-\iop t B} ) \psi\|_2
    \leq |t|\|A-B\| \|\psi\|_2
    \leq |t|\|A-B\| \|\psi\|_1.
\end{equation}

\begin{proof}[Proof of Theorem~\ref{t:PT}]
    Let $H \in \PT$ with $\vertiii{H} =R<\infty$ be given, choose $H_k$ of period $q_k$ satisfying \eqref{eq:PTdef}, fix $\varepsilon>0$, and put
    \begin{equation}
        T_k = e^{\varepsilon^{-1}C(R)q_k}
    \end{equation}
    with $C(R)$ the constant from Theorem~\ref{t:quant:dispersion}. This choice gives the following for $t \geq T_k$:
    \begin{equation} \label{eq:PTtgeqTk}
        \langle t \rangle^{\varepsilon} 
        \geq T_k^{\varepsilon} = e^{C(R)q_k}.
    \end{equation}
    Likewise, for $k$ large enough, \eqref{eq:PTdef} yields
    \begin{equation} \label{eq:PTtleqTk+1}
        \|H - H_k\|
        \leq e^{-\frac{4C(R)}{3\varepsilon} q_{k+1}}
        = T_{k+1}^{-4/3}.
    \end{equation}
    Combining \eqref{eq:duhamel}, \eqref{eq:PTtgeqTk}, and \eqref{eq:PTtleqTk+1}, we obtain the following for $t \in [T_k, T_{k+1}]$ and $k$ large enough:
    \begin{align*}
\|e^{-{\iop}tH}\psi\|_\infty
& \leq \|e^{-{\iop}tH_k}\psi\|_\infty + \|(e^{-{\iop}tH}-e^{-{\iop}tH_k})\psi\|_\infty 
\\
& \leq e^{C(R)q_k} \langle t \rangle^{-1/3}\|\psi\|_1 + T_{k+1} \|H - H_k\| \|\psi\|_1 \\
& \leq \langle t \rangle^{\varepsilon} \langle t \rangle^{-1/3} \|\psi\|_1 + T_{k+1}^{-1/3}\|\psi\|_1 \\
& \lesssim \langle t \rangle^{-\left( \frac13 - \varepsilon \right)} \|\psi \|_1,
\end{align*}
concluding the argument. 
\end{proof}

\begin{remark} \label{rem:EC}
    The argument above also gives a conclusion for Jacobi matrices $H$ with $q_k$-periodic approximants $H_k$ satisfying
    \begin{equation} \label{eq:ECeta}
        \|H - H_k\| \leq Ae^{-\eta q_{k+1}}
    \end{equation}
    for constants $A,\eta>0$, provided $\eta$ is not too small. Indeed, choosing $R$ with $\vertiii{H_k}\leq R$ for all $k$ and supposing \eqref{eq:ECeta} holds for some $\eta > 3C(R)$, reprising the argument above with $T_k = e^{\frac{3}{4}(C(R) + \eta) q_{k}}$ gives
        \begin{align*}
\|e^{-{\iop}tH}\psi\|_\infty
& \leq \|e^{-{\iop}tH_k}\psi\|_\infty + \|(e^{-{\iop}tH}-e^{-{\iop}tH_k})\psi\|_\infty 
\\
& \leq e^{C(R)q_k} \langle t \rangle^{-1/3}\|\psi\|_1 + T_{k+1} \|H - H_k\| \|\psi\|_1 \\
& \leq \langle t \rangle^{\frac{4C(R)}{3(C(R) + \eta)}} \langle t \rangle^{-1/3} \|\psi\|_1 
+ AT_{k+1} T_{k+1}^{-\frac{4\eta}{3(C(R)+\eta)}} \|\psi\|_1 \\
& \lesssim \langle t \rangle^{- \frac{\eta - 3C(R)}{3C(R) + 3\eta}} \|\psi \|_1,
\end{align*}
so one has dispersion at a rate dependent on $\eta$ (and $R$). One also sees that if $\eta$ can be taken arbitrarily large, then the dispersion rate can be made as close to $t^{-1/3}$ as desired.
\end{remark}

The following lemma will be useful in the proof of Theorem~\ref{t:LPrecDisp}. The proof is elementary and included for the reader's convenience.
\begin{lemma} \label{pp:recur}
    If $H: \ell^2(\bbZ) \to \ell^2(\bbZ)$ is a self-adjoint operator with pure point spectrum, then for any $\phi$ with $\|\phi\|=1$, one has
    \begin{equation} \label{eq:pprecurrent}
        \limsup_{t \to \infty} |\langle \phi, e^{-\iop t H} \phi \rangle| = 1.
    \end{equation}
\end{lemma}

\begin{proof}
    Let $\{\varphi_j : j \in \bbN\}$ be a complete orthonormal set of eigenvectors of $H$ with corresponding eigenvalues $E_j$, fix $\varepsilon>0$, and choose $N$ large enough that
    \begin{equation}
        \sum_{j=N+1}^\infty |\langle \varphi_j, \phi\rangle|^2 < \varepsilon/2.
    \end{equation}
    We then have
    \begin{equation}
        |\langle \phi, e^{-\iop t H } \phi \rangle|
        = \left|\sum_{j=1}^\infty e^{-\iop t E_j} |\langle \varphi_j,\phi\rangle|^2 \right|
        \geq \left|\sum_{j=1}^N e^{-\iop t E_j} |\langle \varphi_j,\phi\rangle|^2 \right| - \frac{\varepsilon}{2}.
    \end{equation}
    Choosing $t_k \to \infty$ so that $e^{-\iop t_k E_j} \to 1$ for each $j=1,\ldots,N$ gives
    \begin{align*}
        \limsup_{t \to \infty} |\langle \phi, e^{-\iop t H} \phi \rangle|
        & \geq \limsup_{k \to \infty} \left|\sum_{j=1}^N e^{-\iop t_k E_j} |\langle \varphi_j,\phi\rangle|^2 \right| - \frac{\varepsilon}{2} \\
        & \geq  \left|\sum_{j=1}^N  |\langle \varphi_j,\phi\rangle|^2 \right| - \frac{\varepsilon}{2} \\
        & \geq 1 - \varepsilon.
    \end{align*}
    As this holds for arbitrary $\varepsilon>0$ (and $|\langle \phi, e^{-\iop t H}\phi \rangle| \leq \|\phi\|^2 =1$ for all $t$ by the Cauchy--Schwarz inequality), the proof of \eqref{eq:pprecurrent} is complete.
\end{proof}

\begin{proof}[Proof of Theorem~\ref{t:LPrecDisp}]
Let $\varphi$ be given, and define for each $k \in \bbN$:
\begin{align*}
U(k) = \bigcup_{t\geq k} \{v\in\LP : \|e^{-\iop tH_v}\|_{\ell^1(\bbZ) \to \ell^\infty(\bbZ)} < t^{-1/3} \varphi(t) /k\}.
\end{align*}
 Each $U(k)$ is trivially open and is dense by \cite{DamFilYou2026FMS}, so $\calD  = \calD_\varphi = \bigcap_k U(k)$ is then a dense $G_\delta$ set.
If $v \in \calD$, then there are $t_k \to \infty$ 
satisfying
\begin{equation}
    \|e^{-\iop t_k H_v}\|_{\ell^1(\bbZ) \to \ell^\infty(\bbZ)} 
    < \tfrac1k t_k^{-1/3} \varphi(t_k),
\end{equation}
so \eqref{eq:subseqDispersion}  holds on a dense $G_\delta$ subset of $\LP$.

Next, we show that \eqref{eq:subseqRecurrence} holds generically. For $T,\varepsilon>0$, let 
\begin{equation}
    \calR_{T,\varepsilon} := \bigcup_{t \geq T} \{ v \in \LP(\bbZ)  : |\langle \delta_0, e^{-\iop t H_v} \delta_0 \rangle| > 1-\varepsilon\}.
\end{equation}
On one hand, $\calR_{T, \varepsilon}$ is open.
On the other hand, it is dense because (by Lemma~\ref{pp:recur}) it contains any $v \in \LP$ for which $H_v$ has pure point spectrum, which in turn is dense by \cite{DamGor2016AdvM}.
Intersecting over a sequence $T\to \infty$ and $\varepsilon\downarrow 0$ gives a dense $G_\delta$ set $\calR$ on which \eqref{eq:subseqRecurrence} holds.
Intersecting the two residual sets gives yet another residual set upon which both hold.

Finally, there is another residual set $\calS \subseteq \LP$ with the property that $H_{v}$ exhibits singular continuous spectrum for all $v \in \mathcal{S}$. This is essentially due to Avila \cite{Avi2009CMP}; a proof in this formulation can be found as \cite[Theorem~1.1]{FGH2026JFA} with $k=1$.
Consequently, $\calG:=\calD\cap \calR \cap \calS$ is a residual set (as it is the intersection of three residual sets) such that for every $v \in \calG$, $H_{v}$ has purely singular continuous spectrum and a recurrent dispersive estimate of the required form.
\end{proof}

\begin{appendix}
    \section{Lower Bound on Band Length}
    The purpose of this appendix is to give a simple proof of the lower bound from \eqref{eq:bandLB}.

    \begin{theorem} \label{t:bandLB}
        For each $R>0$, there is a constant $B = B(R)>0$ with the following property.
        For any $q$-periodic Jacobi matrix $H= H_{a,v}$ satisfying $\vertiii{H} \leq R$, one has
        \begin{equation}
            |\sigma_j| \geq \frac{1}{e^{Bq}}.
        \end{equation}
    \end{theorem}

    \begin{proof}
    This can be proved by the method of Last \cite{Last1992CMP} (which gives better constants), but the following simple argument suffices for our purposes.
The assumptions $\|a\|_\infty, \|v\|_\infty \leq R$ imply that  $|E|\leq 3R$ for all $E \in \spectrum H_{a,v}$. Since $\|a^{-1}\|_\infty \leq R$, we have the following for any $n$ and any $E \in \spectrum H$: 
\begin{equation}
    \left\| \frac{1}{a_n} \begin{bmatrix}
        E - v_n & -1 \\  a_n^2 & 0
    \end{bmatrix} \right\|
    \leq 4R^2 + 2R.
\end{equation}
Abbreviate $\onestep_n(E) := \onestep(a_n,v_n,E)$.
By the product rule, we have the following estimate for $E \in \spectrum H_{a,v}$:
\begin{align*}
    &|\Delta'(E)|
     = \left| \trace \monodromy'(E) \right| \\[2mm]
    &\qquad  = \left| \trace \left( \sum_{n=1}^q \onestep_q(E) \cdots \onestep_{n+1}(E) \frac{1}{a_n} \begin{bmatrix}
        1 & 0 \\ 0 & 0
    \end{bmatrix} \onestep_{n-1}(E) \cdots \onestep_1(E) \right) \right| \\[2mm]
    & \qquad \leq qR (4R^2 + 2R)^{q-1}.
\end{align*}
Since $\Delta$ increases from $-2$ to $2$ or decreases from $2$ to $-2$ on each band, we arrive at
        \begin{equation}
            |\sigma_j| \geq
            \frac{4}{qR(4R^2 + 2R)^{q-1}},
        \end{equation}
        so the result follows with an appropriate choice of $B$.
    \end{proof}
\end{appendix}

\section*{Declarations}
\subsection*{Conflict of Interest} The authors have no conflicts of interest to declare.
\subsection*{Funding Information}
D.\ D.\ was supported in part by NSF grant DMS--2349919. J.\ F.\ was supported in part by the National Science Foundation under grant DMS--2513006 and by the  Simons Foundation under grant TSM-00013720.
This material is based upon work supported by the National Science Foundation under Award No. 2503339.
\subsection*{AI Usage}
AI tools, specifically, ChatGPT Astra, were used to obtain the proof of Theorem~\ref{thm:eqConvexity}. The geometric proof presented here represents a distillation of several obscure AI-generated arguments.
All text in this article was written by the authors, who assume all responsibility for correctness of arguments, accuracy of citations, and readability of prose.


\bibliographystyle{amsalpha}

\bibliography{refs}

@article{AmmariSabri2022AMP,
  author = {Ammari, Ka{\"i}s and Sabri, Mostafa},
  title = {Dispersion for {S}chr\"odinger operators on regular trees},
  journal = {Anal. Math. Phys.},
  fjournal = {Analysis and Mathematical Physics},
  volume = {12},
  year = {2022},
  number = {2},
  pages = {Paper No. 56, 38},
  issn = {1664-2368,1664-235X},
  mrclass = {34B45 (46N50 81Q20 81Q35)},
  mrnumber = {4396664},
  mrreviewer = {Juan\ Manuel\ P\'erez-Pardo},
  doi = {10.1007/s13324-022-00664-y},
}

@article{AschKnauf1998Nonlin,
  author = {Joachim Asch and Andreas  Knauf},
  title = {Motion in periodic potentials},
  journal = {Nonlinearity},
  year = {1998},
  fjournal = {Nonlinearity},
  volume = {11},
  number = {1},
  pages = {175--200},
  issn = {0951-7715},
  mrclass = {58F11 (34C11 35Q40 58F05 70H05 81Q20)},
  mrnumber = {1492956},
  mrreviewer = {Gabriel P. Paternain},
  doi = {10.1088/0951-7715/11/1/011},
  url = {https://doi.org/10.1088/0951-7715/11/1/011},
}

@article{Avi2009CMP,
  author = {Artur Avila},
  title = {On the spectrum and {L}yapunov exponent of limit periodic {S}chr\"{o}dinger operators},
  journal = {Comm. Math. Phys.},
  fjournal = {Communications in Mathematical Physics},
  volume = {288},
  year = {2009},
  number = {3},
  pages = {907--918},
  issn = {0010-3616},
  mrclass = {47B36 (37B25 39A70 47B80)},
  mrnumber = {2504859},
  mrreviewer = {David Damanik},
  doi = {10.1007/s00220-008-0667-2},
  url = {https://doi.org/10.1007/s00220-008-0667-2},
}

@article {BambusiZhao2020Advances,
    AUTHOR = {Bambusi, Dario and Zhao, Zhiyan},
     TITLE = {Dispersive estimate for quasi-periodic {S}chr\"odinger
              operators on 1-{$d$} lattices},
   JOURNAL = {Adv. Math.},
  FJOURNAL = {Advances in Mathematics},
    VOLUME = {366},
      YEAR = {2020},
     PAGES = {107071, 27},
      ISSN = {0001-8708,1090-2082},
   MRCLASS = {35Q55 (35R02)},
  MRNUMBER = {4070305},
MRREVIEWER = {Jing\ Bo\ Chen},
       DOI = {10.1016/j.aim.2020.107071},
       URL = {https://doi-org.srv-proxy1.library.tamu.edu/10.1016/j.aim.2020.107071},
}

@article {BenkoDragnev2012ConstrApp,
    AUTHOR = {Benko, David and Dragnev, Peter D.},
     TITLE = {Balayage ping-pong: a convexity of equilibrium measures},
   JOURNAL = {Constr. Approx.},
  FJOURNAL = {Constructive Approximation. An International Journal for
              Approximations and Expansions},
    VOLUME = {36},
      YEAR = {2012},
    NUMBER = {2},
     PAGES = {191--214},
      ISSN = {0176-4276,1432-0940},
   MRCLASS = {31A15 (33C45 41A10 42C05)},
  MRNUMBER = {2957308},
MRREVIEWER = {Dmitriy\ Karp},
       DOI = {10.1007/s00365-011-9143-x},
       }

@article {BenDraTot2012RevMatIbe,
    AUTHOR = {Benko, David and Dragnev, Peter and Totik, Vilmos},
     TITLE = {Convexity of harmonic densities},
   JOURNAL = {Rev. Mat. Iberoam.},
  FJOURNAL = {Revista Matem\'atica Iberoamericana},
    VOLUME = {28},
      YEAR = {2012},
    NUMBER = {4},
     PAGES = {947--960},
      ISSN = {0213-2230,2235-0616},
   MRCLASS = {31A15},
  MRNUMBER = {2990129},
MRREVIEWER = {Mehrdad\ Simkani},
       DOI = {10.4171/RMI/698},
}

@incollection{BoutetdeMonvelSabri2022OTAA,
  author = {Boutet de Monvel, Anne and Sabri, Mostafa},
  title = {Ballistic transport in periodic and random media},
  booktitle = {From complex analysis to operator theory---a panorama},
  publisher = {Birkh\"auser/Springer, Cham},
  year = {2023},
  series = {Oper. Theory Adv. Appl.},
  volume = {291},
  pages = {163--216},
  isbn = {978-3-031-31138-3; 978-3-031-31139-0},
  mrclass = {47A10 (35J10 47A55 81Q10)},
  mrnumber = {4651275},
  doi = {10.1007/978-3-031-31139-0\_10},
  url = {https://doi.org/10.1007/978-3-031-31139-0_10},
}

@article {Bucaj2018OperMat,
    AUTHOR = {Bucaj, Valmir},
     TITLE = {The {K}unz-{S}ouillard approach to localization for {J}acobi
              operators},
   JOURNAL = {Oper. Matrices},
  FJOURNAL = {Operators and Matrices},
    VOLUME = {12},
      YEAR = {2018},
    NUMBER = {4},
     PAGES = {1099--1127},
      ISSN = {1846-3886,1848-9974},
   MRCLASS = {47B37},
  MRNUMBER = {3890753},
MRREVIEWER = {Shiv\ Kumar\ Kaushik},
       DOI = {10.7153/oam-2018-12-67},
       URL = {https://doi.org/10.7153/oam-2018-12-67},
}

@article{Cai2006DPDE,
  author = {Cai, Kaihua},
  title = {Dispersion for {S}chr\"odinger operators with one-gap periodic potentials on {$\mathbb{R}^1$}},
  journal = {Dyn. Partial Differ. Equ.},
  fjournal = {Dynamics of Partial Differential Equations},
  volume = {3},
  year = {2006},
  number = {1},
  pages = {71--92},
  issn = {1548-159X,2163-7873},
  mrclass = {35Q40 (35B10 37K15 47B80 47E05)},
  mrnumber = {2221747},
  mrreviewer = {Fran\c cois\ Germinet},
  doi = {10.4310/DPDE.2006.v3.n1.a2},
  url = {https://doi.org/10.4310/DPDE.2006.v3.n1.a2},
  eprint = {arXiv:math-ph/0411002},
}

@article{cataliotti2001josephson,
  title={Josephson junction arrays with {B}ose-{E}instein condensates},
  author={Cataliotti, Francesco Saverio and Burger, Sven and Fort, Chiara and Maddaloni, Pasquale  and Minardi, Francesco and Trombettoni, Andrea and Smerzi, Augusto and Inguscio, Massimo},
  journal={Science},
  volume={293},
  number={5531},
  pages={843--846},
  year={2001},
  publisher={American Association for the Advancement of Science}
}

@article{cataliotti2003superfluid,
  title={Superfluid current disruption in a chain of weakly coupled {B}ose--{E}instein condensates},
  author={Cataliotti, Francesco Saverio and Fallani, Leonardo and Ferlaino, Francesca  and Fort, Chiara and Maddaloni, Pasquale  and Inguscio, Massimo},
  journal={New J. Physics},
  volume={5},
  number={1},
  pages={71},
  year={2003},
  publisher={IOP Publishing}
}

@article{Cuccagna2008CPDE,
  author = {Cuccagna, Scipio},
  title = {Dispersion for {S}chr{\"o}dinger equation with periodic potential in 1D},
  journal = {Commun. Partial. Differ. Equ.},
  year = {2008},
  volume = {33},
  number = {11},
  pages = {2064--2095},
  publisher = {Taylor \& Francis},
}

@article{Cuccagna2009SIAMJMA,
  author = {Cuccagna, Scipio and Tarulli, Mirko},
  title = {On asymptotic stability of standing waves of discrete {S}chr{\"o}dinger equation in {$\mathbb{Z}$} },
  journal = {SIAM J. Math. Anal.},
  year = {2009},
  volume = {41},
  number = {3},
  pages = {861--885},
  publisher = {SIAM},
}

@incollection{DamFil2020LPSurvey,
  author = {Damanik, David and Fillman, Jake},
  title = {Spectral properties of limit-periodic operators},
  booktitle = {Analysis and geometry on graphs and manifolds},
  series = {London Math. Soc. Lecture Note Ser.},
  volume = {461},
  pages = {382--444},
  publisher = {Cambridge Univ. Press, Cambridge},
  year = {2020},
  isbn = {978-1-108-71318-4},
  mrclass = {47B36 (39A70)},
  mrnumber = {4412981},
}

@book{DamFil2022ESO1,
  author = {Damanik, David and Fillman, Jake},
  title = {One-dimensional ergodic {S}chr\"odinger operators---{I}. {G}eneral theory},
  series = {Graduate Studies in Mathematics},
  volume = {221},
  publisher = {American Mathematical Society, Providence, RI},
  year = {2022},
  pages = {xv+444},
  isbn = {978-1-4704-5606-1; [9781470470869]; [9781470470852]},
  mrclass = {47-01 (28Axx 35Q41 37Axx 47B36 81Q10 82B44)},
  mrnumber = {4567742},
  mrreviewer = {Christian\ Seifert},
}

@book{DamFil2024ESO2,
  author = {Damanik, David and Fillman, Jake},
  title = {One-dimensional ergodic {S}chr\"odinger operators---{II}. {S}pecific classes},
  series = {Graduate Studies in Mathematics},
  volume = {249},
  publisher = {American Mathematical Society, Providence, RI},
  year = {2024},
  pages = {vii--xvi and 445--1068},
  isbn = {978-1-4704-6503-2},
  mrclass = {47-01 (35J10 47B36 81Q10 82B44)},
  mrnumber = {4840232},
  mrreviewer = {Christian\ Seifert},
}

@article{DamFilLuk2017JST,
  author = {Damanik, David and Fillman, Jake and Lukic, Milivoje},
  title = {Limit-periodic continuum {S}chr\"{o}dinger operators with zero measure {C}antor spectrum},
  journal = {J. Spectr. Theory},
  fjournal = {Journal of Spectral Theory},
  volume = {7},
  year = {2017},
  number = {4},
  pages = {1101--1118},
  issn = {1664-039X},
  mrclass = {34L40 (35J10 47B36)},
  mrnumber = {3737889},
  mrreviewer = {Rafael R. Del Rio},
  doi = {10.4171/JST/186},
}

@article{DamFilWan2023MathNach,
  author = {Damanik, David and Fillman, Jake and Wang, Chunyi},
  title = {Thin spectra and singular continuous spectral measures for limit-periodic {J}acobi matrices},
  journal = {Math. Nachr.},
  fjournal = {Mathematische Nachrichten},
  volume = {296},
  year = {2023},
  number = {9},
  pages = {4279--4297},
  issn = {0025-584X,1522-2616},
  mrclass = {47B36 (15B05 35J10 47A10 58J51)},
  mrnumber = {4644562},
  mrreviewer = {Dolores\ Barrios Rolan\'ia},
  doi = {10.1002/mana.202100561},
  url = {https://doi-org.srv-proxy1.library.tamu.edu/10.1002/mana.202100561},
}

@article{DamFilYou2026FMS,
    author = {David Damanik and Jake Fillman and Giorgio Young},
   title={Optimal dispersion for discrete periodic {S}chr\"odinger operators}, 
      year={2026},
      note={arxiv:2505.14475},
      journal = {Forum Math. Sigma},
      pages = {1-14},
}

@article{DamGor2016AdvM,
  author = {David Damanik and Anton Gorodetski},
  title = {An extension of the {K}unz-{S}ouillard approach to localization in one dimension and applications to almost-periodic {S}chr\"{o}dinger operators},
  journal = {Adv. Math.},
  fjournal = {Advances in Mathematics},
  volume = {297},
  year = {2016},
  pages = {149--173},
  issn = {0001-8708},
  mrclass = {47A35},
  mrnumber = {3498796},
  mrreviewer = {Aleksey Kostenko},
  doi = {10.1016/j.aim.2016.04.006},
  url = {https://doi.org/10.1016/j.aim.2016.04.006},
}

@article{DamLukYes2015CMP,
  author = {David Damanik and Milivoje Luki\'c and William Yessen},
  title = {Quantum dynamics of periodic and limit-periodic {J}acobi and block {J}acobi matrices with applications to some quantum many body problems},
  journal = {Comm. Math. Phys.},
  fjournal = {Communications in Mathematical Physics},
  volume = {337},
  year = {2015},
  number = {3},
  pages = {1535--1561},
  issn = {0010-3616},
  mrclass = {81U10},
  mrnumber = {3339185},
  mrreviewer = {Marco Falconi},
  doi = {10.1007/s00220-015-2346-4},
  url = {https://doi.org/10.1007/s00220-015-2346-4},
}

@article {DamMalYou2026JST,
    AUTHOR = {Damanik, David and Malinovitch, Tal and Young, Giorgio},
     TITLE = {What is ballistic transport?},
   JOURNAL = {J. Spectr. Theory},
  FJOURNAL = {Journal of Spectral Theory},
    VOLUME = {16},
      YEAR = {2026},
    NUMBER = {3},
     PAGES = {995--1032},
      ISSN = {1664-039X,1664-0403},
   MRCLASS = {35Q41 (35J10 35Q40 81Q10)},
  MRNUMBER = {5112773},
       DOI = {10.4171/jst/532},
       URL = {https://doi.org/10.4171/jst/532},
}

@article{efremidis2002discrete,
  title={Discrete solitons in photorefractive optically induced photonic lattices},
  author={Efremidis, Nikos K. and Sears, Suzanne and Christodoulides, Demetrios N. and Fleischer, Jason W. and Segev, Mordechai},
  journal={Phys. Rev. E},
  volume={66},
  number={4},
  pages={046602},
  year={2002},
  publisher={APS}
}

@article {Egorova1987,
    AUTHOR = {Egorova, Iryna E.},
     TITLE = {Spectral analysis of {J}acobi limit-periodic matrices},
   JOURNAL = {Dokl. Akad. Nauk Ukrain. SSR Ser. A},
  FJOURNAL = {Akademiya Nauk Ukrainsko\u i\ SSR. Doklady. Seriya A.
              Fiziko-Matematicheskie i Tekhnicheskie Nauki},
      YEAR = {1987},
    NUMBER = {3},
     PAGES = {7--9, 85},
      ISSN = {0201-8446},
   MRCLASS = {15A18 (15A60)},
  MRNUMBER = {892383},
MRREVIEWER = {Lajos\ L\'aszl\'o},
}

@article{EgoKopTes2015JST,
  title={Dispersion estimates for one-dimensional discrete {S}chr{\"o}dinger and wave equations},
  author={Egorova, Iryna and Kopylova, Elena  A. and Teschl, Gerald},
  journal={J.  Spectr. Theory},
  volume={5},
  number={4},
  pages={663--696},
  year={2015},
}

@article{EFGL2022JFA,
  author = {Eichinger, Benjamin and Fillman, Jake and Gwaltney, Ethan and Luki\'c, Milivoje},
  title = {Limit-periodic {D}irac operators with thin spectra},
  journal = {J. Funct. Anal.},
  fjournal = {Journal of Functional Analysis},
  volume = {283},
  year = {2022},
  number = {12},
  pages = {Paper No. 109711, 35},
  issn = {0022-1236,1096-0783},
  mrclass = {34L16 (47A10 47B36)},
  mrnumber = {4489275},
  doi = {10.1016/j.jfa.2022.109711},
  url = {https://doi.org/10.1016/j.jfa.2022.109711},
}

@article{eisenberg2002optical,
  title={Optical discrete solitons in waveguide arrays. I. Soliton formation},
  author={Eisenberg, Hagai S. and Morandotti, R. and Silberberg, Yaron and Arnold, J.M. and Pennelli, Giovanni and Aitchison, J.S.},
  journal={J. Opt. Soc. Am. B},
  volume={19},
  number={12},
  pages={2938--2944},
  year={2002},
  publisher={Optical Society of America}
}

@article{Fillman2017CMP,
  author = {Jake Fillman},
  title = {Ballistic transport for limit-periodic {J}acobi matrices with applications to quantum many-body problems},
  journal = {Comm. Math. Phys.},
  fjournal = {Communications in Mathematical Physics},
  volume = {350},
  year = {2017},
  number = {3},
  pages = {1275--1297},
  issn = {0010-3616},
  mrclass = {81Q10},
  mrnumber = {3607475},
  doi = {10.1007/s00220-016-2785-6},
}

@incollection{Fillman2021OTAA,
  author = {Fillman, Jake},
  title = {Ballistic transport for periodic {J}acobi operators on {$\mathbb{Z}^d$}},
  booktitle = {From operator theory to orthogonal polynomials, combinatorics, and number theory---a volume in honor of {L}ance {L}ittlejohn's 70th birthday},
  publisher = {Birkh\"{a}user/Springer, Cham},
  year = {2021},
  series = {Oper. Theory Adv. Appl.},
  volume = {285},
  pages = {57--68},
  mrclass = {47B36 (82C10)},
  mrnumber = {4367462},
}

@article{FGH2026JFA,
  author = {Fillman, Jake and Gradner, Michala N. and Hendricks, Hannah J.},
  title = {Thin spectra for periodic and ergodic word models},
  journal = {J. Funct. Anal.},
  fjournal = {Journal of Functional Analysis},
  volume = {290},
  year = {2026},
  number = {9},
  pages = {Paper No. 111385, 24},
  issn = {0022-1236,1096-0783},
  mrclass = {47B36 (37A25 47A10)},
  mrnumber = {5030423},
  mrreviewer = {Kenan\ Ta\c s},
  doi = {10.1016/j.jfa.2026.111385},
}

@article {FillmanLukic2017JST,
    AUTHOR = {Fillman, Jake and Lukic, Milivoje},
     TITLE = {Spectral homogeneity of limit-periodic {S}chr\"odinger
              operators},
   JOURNAL = {J. Spectr. Theory},
  FJOURNAL = {Journal of Spectral Theory},
    VOLUME = {7},
      YEAR = {2017},
    NUMBER = {2},
     PAGES = {387--406},
      ISSN = {1664-039X,1664-0403},
   MRCLASS = {34L40 (35J10 47B36)},
  MRNUMBER = {3662012},
MRREVIEWER = {Maxim\ Zinchenko},
       DOI = {10.4171/JST/166},
       URL = {https://doi.org/10.4171/JST/166},
}

@article{FillmanOng2017JFA,
  author = {Fillman, Jake and Ong, Darren C.},
  title = {Purely singular continuous spectrum for limit-periodic {CMV} operators with applications to quantum walks},
  journal = {J. Funct. Anal.},
  fjournal = {Journal of Functional Analysis},
  volume = {272},
  year = {2017},
  number = {12},
  pages = {5107--5143},
  issn = {0022-1236,1096-0783},
  mrclass = {47L10},
  mrnumber = {3639523},
  mrreviewer = {Miguel\ Lacruz},
  doi = {10.1016/j.jfa.2017.01.021},
  url = {https://doi.org/10.1016/j.jfa.2017.01.021},
}

@article{Firsova1987MatSb,
  title = {The direct and inverse scattering problems for the one-dimensional perturbed {H}ill operator},
  author = {Firsova, Nataliya E.},
  journal = {Mat. Sb.},
  volume = {58},
  number = {2},
  pages = {351},
  year = {1987},
  publisher = {IOP Publishing},
}

@article{Firsova1996JMP,
  author = {Firsova, Nataliya  E.},
  title = {On the time decay of a wave packet in a one-dimensional finite band periodic lattice},
  journal = {J. Math. Phys.},
  fjournal = {Journal of Mathematical Physics},
  volume = {37},
  year = {1996},
  number = {3},
  pages = {1171--1181},
  issn = {0022-2488,1089-7658},
  mrclass = {81Q05 (34L25)},
  mrnumber = {1377626},
  mrreviewer = {Luciano\ Canton},
  doi = {10.1063/1.531454},
  url = {https://doi.org/10.1063/1.531454},
}

@article{GeKachkovskiy2023CPAM,
  author = {Ge, Lingrui and Kachkovskiy, Ilya},
  title = {Ballistic transport for one-dimensional quasiperiodic {S}chr\"{o}dinger operators},
  journal = {Comm. Pure Appl. Math.},
  fjournal = {Communications on Pure and Applied Mathematics},
  volume = {76},
  year = {2023},
  number = {10},
  pages = {2577--2612},
  issn = {0010-3640},
  mrclass = {47A35 (35Q55 47B36 81Q10)},
  mrnumber = {4630598},
}

@misc{solis2026noteacspectrumlimitperiodic,
      title={A note on the A.C. spectrum of limit-periodic {J}acobi matrices}, 
      author={Moises Gomez Solis},
      year={2026},
      note={arxiv:2610.19416},
      archivePrefix={arXiv},
      primaryClass={math.SP},
      url={https://arxiv.org/abs/2610.19416}, 
}

@book{GreeneKrantz2002GSM,
  author = {Greene, Robert E. and Krantz, Steven G.},
  title = {Function theory of one complex variable},
  series = {Graduate Studies in Mathematics},
  volume = {40},
  edition = {Second},
  publisher = {American Mathematical Society, Providence, RI},
  year = {2002},
  pages = {xx+502},
  isbn = {0-8218-2905-X},
  mrclass = {30-01},
  mrnumber = {1871333},
  doi = {10.1090/gsm/040},
  url = {https://doi.org/10.1090/gsm/040},
}

@misc{han2025ballistictransportlimitperiodicjacobi,
      title={On the Ballistic Transport for limit-periodic {J}acobi Matrices}, 
      author={Rui Han and Moises Gomez Solis},
      year={2025},
      eprint={2509.16787},
      archivePrefix={arXiv},
      primaryClass={math-ph},
      url={https://arxiv.org/abs/2509.16787}, 
}

@article{kevrekidis2001discrete,
  title={The discrete nonlinear {S}chr{\"o}dinger equation: a survey of recent results},
  author={Kevrekidis, Panayotis G. and Rasmussen, Kim {\O}. and Bishop, Alan R.},
  journal={International Journal of Modern Physics B},
  volume={15},
  number={21},
  pages={2833--2900},
  year={2001},
  publisher={World Scientific}
}

@article{KPST2025CMP,
  author = {Kerner, Joachim and Post, Olaf and Sabri, Mostafa and T\"aufer, Matthias},
  title = {The curious spectra and dynamics of non-locally finite crystals},
  journal = {Comm. Math. Phys.},
  fjournal = {Communications in Mathematical Physics},
  volume = {406},
  year = {2025},
  number = {7},
  pages = {Paper No. 169, 48},
  issn = {0010-3616,1432-0916},
  mrclass = {81Q10 (05E10 37E25 81Q35)},
  mrnumber = {4922843},
  mrreviewer = {Natalia\ Saburova},
  doi = {10.1007/s00220-025-05346-x},
}

@article{KunzSouillard1980CMP,
  author = {Kunz, Herv\'{e} and Souillard, Bernard},
  title = {Sur le spectre des op\'{e}rateurs aux diff\'{e}rences finies al\'{e}atoires},
  journal = {Comm. Math. Phys.},
  fjournal = {Communications in Mathematical Physics},
  volume = {78},
  year = {1980/81},
  number = {2},
  pages = {201--246},
  issn = {0010-3616,1432-0916},
  mrclass = {39A12 (60H25 82A05)},
  mrnumber = {597748},
  url = {http://projecteuclid.org/euclid.cmp/1103908590},
}

@article{Last1992CMP,
  author = {Last, Yoram},
  title = {On the measure of gaps and spectra for discrete {$1$}{D} {S}chr\"{o}dinger operators},
  journal = {Comm. Math. Phys.},
  fjournal = {Communications in Mathematical Physics},
  volume = {149},
  year = {1992},
  number = {2},
  pages = {347--360},
  issn = {0010-3616},
  mrclass = {47B39 (47N50 81Q10 82B10)},
  mrnumber = {1186033},
  mrreviewer = {Russell A. Johnson},
  url = {http://projecteuclid.org/euclid.cmp/1104251226},
}

@misc{leclerc2026fourierdecayequilibriumstates,
      title={Fourier decay of equilibrium states and the {F}ibonacci {H}amiltonian}, 
      author={Gaétan Leclerc},
      year={2026},
      note={arxiv:2507.23731},
      archivePrefix={arXiv},
      primaryClass={math.DS},
      url={https://arxiv.org/abs/2507.23731}, 
}

@misc{LeclercSabriSahlstenPreprint,
      title={Arbitrarily Fast Quantum Dispersion in Long-Range Crystals}, 
      author={Ga\'etan Leclerc and Mostafa Sabri and Tuomas Sahlsten},
      year={2026},
      note={arxiv:2608.27326},
      archivePrefix={arXiv},
      primaryClass={math.SP},
      url={https://arxiv.org/abs/2608.27326}, 
}

@book{Lukic2022book,
  author = {Luki\'c, Milivoje},
  title = {A first course in spectral theory},
  series = {Graduate Studies in Mathematics},
  volume = {226},
  publisher = {American Mathematical Society, Providence, RI},
  year = {2022},
  pages = {xvi+472},
  isbn = {[9781470466565]; [9781470471927]; [9781470471910]},
  mrclass = {47-01 (34B07 47B15 47B25 47B36)},
  mrnumber = {4497147},
  mrreviewer = {Valentin\ Keyantuo},
  doi = {10.1090/gsm/226},
  url = {https://doi.org/10.1090/gsm/226},
}

@article {MiZhao2020JMAA,
    AUTHOR = {Mi, Yue and Zhao, Zhiyan},
     TITLE = {Dispersive estimate for two-periodic discrete one-dimensional
              {S}chr\"odinger operator},
   JOURNAL = {J. Math. Anal. Appl.},
  FJOURNAL = {Journal of Mathematical Analysis and Applications},
    VOLUME = {485},
      YEAR = {2020},
    NUMBER = {1},
     PAGES = {123768, 15},
      ISSN = {0022-247X,1096-0813},
   MRCLASS = {47D08 (81Q80)},
  MRNUMBER = {4042637},
MRREVIEWER = {Jean-Claude\ Cuenin},
       DOI = {10.1016/j.jmaa.2019.123768},
       URL = {https://doi-org.srv-proxy1.library.tamu.edu/10.1016/j.jmaa.2019.123768},
}

@article {MiZhao2022PAMS,
    AUTHOR = {Mi, Yue and Zhao, Zhiyan},
     TITLE = {Dispersive estimates for periodic discrete one-dimensional
              {S}chr\"odinger operators},
   JOURNAL = {Proc. Amer. Math. Soc.},
  FJOURNAL = {Proceedings of the American Mathematical Society},
    VOLUME = {150},
      YEAR = {2022},
    NUMBER = {1},
     PAGES = {267--277},
      ISSN = {0002-9939,1088-6826},
   MRCLASS = {35P25 (47B36)},
  MRNUMBER = {4335875},
MRREVIEWER = {Petru\ A.\ Cojuhari},
       DOI = {10.1090/proc/15699}, 
}

@article {PasturTkachenko1,
	AUTHOR = {Pastur, Leonid A. and Tkachenko, Vadim A.},
	TITLE = {On the Spectral Theory of the One-Dimensional {S}chr\"odinger Operator with Limit-Periodic Potential ({R}ussian).},
	JOURNAL = {Dokl. Akad. Nauk SSSR},
	VOLUME = {279},
	YEAR = {1984},
	PAGES = {1050–1053},
	MRCLASS = {},
	MRNUMBER = {},
	MRREVIEWER = {},
	DOI = {},
	URL = {},
}

@article {PasturTkachenko2,
	AUTHOR = {Pastur, Leonid A. and Tkachenko, Vadim A.},
	TITLE = {Spectral theory of a class of one-dimensional {S}chr\"odinger operators with limit-periodic potentials ({R}ussian).},
	JOURNAL = { Trudy Moskov. Mat. Obshch.},
	VOLUME = {51},
	YEAR = {1988},
	PAGES = {114–168},
	MRCLASS = {},
	MRNUMBER = {},
	MRREVIEWER = {},
	DOI = {},
	URL = {},
}

@article{peschel2002optical,
  title={Optical discrete solitons in waveguide arrays. 2. Dynamic properties},
  author={Peschel, Ulf and Morandotti, Roberto and Arnold, John M. and Aitchison, J. Stewart and Eisenberg, Hagai S. and Silberberg, Yaron and Pertsch, Thomas and Lederer, Falk},
  journal={J. Opt. Soc. Am. B },
  volume={19},
  number={11},
  pages={2637--2644},
  year={2002},
  publisher={Optical Society of America}
}

@article {Saff2010SurvApproxTh,
    AUTHOR = {Saff, Edward B.},
     TITLE = {Logarithmic potential theory with applications to
              approximation theory},
   JOURNAL = {Surv. Approx. Theory},
  FJOURNAL = {Surveys in Approximation Theory},
    VOLUME = {5},
      YEAR = {2010},
     PAGES = {165--200},
      ISSN = {1555-578X},
   MRCLASS = {31-02 (30E10 31A15 41A20)},
  MRNUMBER = {2734198},
MRREVIEWER = {Anders\ Gustafsson},
}

@book {SaffTotik2017,
    AUTHOR = {Saff, Edward B. and Totik, Vilmos},
     TITLE = {Logarithmic potentials with external fields},
    SERIES = {Grundlehren der mathematischen Wissenschaften},
    VOLUME = {316},
   EDITION = {Second},
      NOTE = {With an appendix by Thomas Bloom},
 PUBLISHER = {Springer, Cham},
      YEAR = {2024},
     PAGES = {xvi+594},
      ISBN = {978-3-031-65132-8; 978-3-031-65133-5},
   MRCLASS = {31-02 (31Axx 32U05)},
  MRNUMBER = {4807484},
       DOI = {10.1007/978-3-031-65133-5},
}

@article {SKW2026JMAA,
    AUTHOR = {Sagiv, Amir and Kassem, Remy and Weinstein, Michael I.},
     TITLE = {Dispersive decay estimates for periodic {J}acobi operators on
              the half-line},
   JOURNAL = {J. Math. Anal. Appl.},
  FJOURNAL = {Journal of Mathematical Analysis and Applications},
    VOLUME = {553},
      YEAR = {2026},
    NUMBER = {1},
     PAGES = {Paper No. 129945, 17},
      ISSN = {0022-247X,1096-0813},
   MRCLASS = {47B36 (35B10 35Q41 42C05 81Q10)},
  MRNUMBER = {4944037},
MRREVIEWER = {Grzegorz\ \'Swiderski},
       DOI = {10.1016/j.jmaa.2025.129945},
}

@incollection {Schlag2007Survey,
    AUTHOR = {Schlag, Wilhelm},
     TITLE = {Dispersive estimates for {S}chr\"odinger operators: a survey},
 BOOKTITLE = {Mathematical aspects of nonlinear dispersive equations},
    SERIES = {Ann. of Math. Stud.},
    VOLUME = {163},
     PAGES = {255--285},
 PUBLISHER = {Princeton Univ. Press, Princeton, NJ},
      YEAR = {2007},
      ISBN = {978-0-691-12955-6; 0-691-12955-X},
   MRCLASS = {35J10 (35-02 35B45 47F05)},
  MRNUMBER = {2333215},
}

@article {Schlag2021JMP,
    AUTHOR = {Schlag, Wilhelm},
     TITLE = {On pointwise decay of waves},
   JOURNAL = {J. Math. Phys.},
  FJOURNAL = {Journal of Mathematical Physics},
    VOLUME = {62},
      YEAR = {2021},
    NUMBER = {6},
     PAGES = {Paper No. 061509, 27},
      ISSN = {0022-2488,1089-7658},
   MRCLASS = {35J10 (47D08 81T20)},
  MRNUMBER = {4276283},
       DOI = {10.1063/5.0042767},
       URL = {https://doi.org/10.1063/5.0042767},
}

@article{Simon1990CMP,
title = {Absence of Ballistic Motion},
journal = {Comm. Math. Phys.},
volume = {149},
number = {},
pages = {209-212},
year = {1990},
issn = {},
doi = {},
url = {},
author = {Barry Simon},
}

@article {Simon2007IPI,
    AUTHOR = {Barry Simon},
     TITLE = {Equilibrium measures and capacities in spectral theory},
   JOURNAL = {Inverse Probl. Imaging},
  FJOURNAL = {Inverse Problems and Imaging},
    VOLUME = {1},
      YEAR = {2007},
    NUMBER = {4},
     PAGES = {713--772},
      ISSN = {1930-8337,1930-8345},
   MRCLASS = {31A15 (34A55 34L05 34L40 35J10 42C05)},
  MRNUMBER = {2350223},
MRREVIEWER = {Franck\ Wielonsky},
       DOI = {10.3934/ipi.2007.1.713},
}

@book{Simon2011:Szego,
  author = {Barry Simon},
  title = {Szeg{\H{o}}'s theorem and its descendants: {S}pectral theory for {$L^2$} perturbations of orthogonal polynomials},
  series = {M. B. Porter Lectures},
  publisher = {Princeton University Press, Princeton, NJ},
  year = {2011},
  pages = {xii+650},
  isbn = {978-0-691-14704-8},
  mrclass = {47B36 (37K15 37K60 42-02 42C05 47A10 47B35)},
  mrnumber = {2743058},
  mrreviewer = {Harry Dym},
}

@article{solis2026two,
  title={Two-scale estimates for limit-periodic Schr\"{o}dinger operators},
  author={Solis, Moises Gomez},
  journal={arXiv preprint arXiv:2611.01303},
  year={2026}
}

@article{sukhorukov2003spatial,
  title={Spatial optical solitons in waveguide arrays},
  author={Sukhorukov, Andrey A. and Kivshar, Yuri S. and Eisenberg, Hagai S. and Silberberg, Yaron},
  journal={IEEE J. Quantum Electron.},
  volume={39},
  number={1},
  pages={31--50},
  year={2003},
  publisher={IEEE}
}

@book{Teschl2009:SpTh,
  author = {Gerald Teschl},
  title = {Mathematical methods in quantum mechanics},
  series = {Graduate Studies in Mathematics},
  volume = {99},
  note = {With applications to Schr\"{o}dinger operators},
  publisher = {American Mathematical Society, Providence, RI},
  year = {2009},
  pages = {xiv+305},
  isbn = {978-0-8218-4660-5},
  mrclass = {81-01 (47N50 81Qxx)},
  mrnumber = {2499016},
  mrreviewer = {Rupert L. Frank},
  doi = {10.1090/gsm/099},
  url = {https://doi.org/10.1090/gsm/099},
}

@article{Young2023JST,
  author = {Young, Giorgio},
  title = {Ballistic transport for limit-periodic {S}chr\"{o}dinger operators in one dimension},
  journal = {J. Spectr. Theory},
  fjournal = {Journal of Spectral Theory},
  volume = {13},
  year = {2023},
  number = {2},
  pages = {451--489},
  issn = {1664-039X,1664-0403},
  mrclass = {47A10 (35J10)},
  mrnumber = {4651730},
  doi = {10.4171/jst/463},
  url = {https://doi.org/10.4171/jst/463},
}

@article{ZhangZhao2017CMP,
  author = {Zhang, Zhiyuan and Zhao, Zhiyan},
  title = {Ballistic transport and absolute continuity of one-frequency {S}chr\"{o}dinger operators},
  journal = {Comm. Math. Phys.},
  fjournal = {Communications in Mathematical Physics},
  volume = {351},
  year = {2017},
  number = {3},
  pages = {877--921},
  issn = {0010-3616},
  mrclass = {35Q41 (35B65 37K55)},
  mrnumber = {3623240},
  mrreviewer = {William R. Green},
  doi = {10.1007/s00220-017-2848-3},
}

@article{Zhao2016CMP,
  author = {Zhao, Zhiyan},
  title = {Ballistic motion in one-dimensional quasi-periodic discrete {S}chr\"{o}dinger equation},
  journal = {Comm. Math. Phys.},
  fjournal = {Communications in Mathematical Physics},
  volume = {347},
  year = {2016},
  number = {2},
  pages = {511--549},
  issn = {0010-3616},
  mrclass = {34L40 (37J40)},
  mrnumber = {3545515},
  mrreviewer = {C\'{e}sar R. de Oliveira},
  doi = {10.1007/s00220-016-2605-z},
}

@book{Stein1993HarmonicAnalysis,
  author    = {Stein, Elias M.},
  title     = {Harmonic Analysis: Real-Variable Methods,
               Orthogonality, and Oscillatory Integrals},
  series    = {Princeton Mathematical Series},
  volume    = {43},
  publisher = {Princeton University Press},
  address   = {Princeton, NJ},
  year      = {1993},
  isbn      = {978-0-691-03216-0}
}

@incollection{BianchiniSalani2013,
  author    = {Bianchini, Massimiliano and Salani, Paolo},
  title     = {Power Concavity for Solutions of Nonlinear
               Elliptic Problems in Convex Domains},
  booktitle = {Geometric Properties for Parabolic and Elliptic PDE's},
  editor    = {Magnanini, Rolando and Sakaguchi, Shigeru
               and Alvino, Angelo},
  series    = {Springer INdAM Series},
  volume    = {2},
  pages     = {35--48},
  year      = {2013},
  publisher = {Springer},
  doi       = {10.1007/978-88-470-2841-8_3}
}

@article{Kennington1985,
  author  = {Kennington, Alan U.},
  title   = {Power Concavity and Boundary Value Problems},
  journal = {Indiana University Mathematics Journal},
  volume  = {34},
  number  = {3},
  pages   = {687--704},
  year    = {1985},
  doi     = {10.1512/iumj.1985.34.34036},
  url     = {https://math.jhu.edu/~js/Math745/kennington.pdf}
}

\end{document}